\documentclass[10pt, a4wide]{article}
\usepackage{graphicx}
\usepackage{amsmath,amsthm,amssymb,enumerate}
\usepackage{color}
\usepackage[utf8]{inputenc}
\catcode`\@=11 \@addtoreset{equation}{section}

\catcode`\@=12
\usepackage{hyperref}
\usepackage{authblk}
\usepackage{etoolbox}
\usepackage{lmodern}
\makeatletter
\patchcmd{\@maketitle}{\LARGE \@title}{\fontsize{18}{6}\selectfont\@title}{}{}
\makeatother

\newtheorem{Theorem}{Theorem}[section]
\newtheorem{Proposition}{Proposition}[section]
\newtheorem{Lemma}{Lemma}[section]
\newtheorem{Corollary}{Corollary}[section]
\newtheorem{Remark}{Remark}[section]
\newtheorem{Definition}{Definition}[section]

\newcommand{\bTheorem}[1]{
	\begin{Theorem} \label{T#1} }
	\newcommand{\eT}{\end{Theorem}}

\newcommand{\bProposition}[1]{
	\begin{Proposition} \label{P#1}}
	\newcommand{\eP}{\end{Proposition}}

\newcommand{\bLemma}[1]{
	\begin{Lemma} \label{L#1} }
	\newcommand{\eL}{\end{Lemma}}

\newcommand{\bCorollary}[1]{
	\begin{Corollary} \label{C#1} }
	\newcommand{\eC}{\end{Corollary}}

\newcommand{\bRemark}[1]{
	\begin{Remark} \label{R#1} }
	\newcommand{\eR}{\end{Remark}}

\newcommand{\bDefinition}[1]{
	\begin{Definition} \label{D#1} }
	\newcommand{\eD}{\end{Definition}}
\newcommand{\fr}{\mathfrak{L}_{\lambda}}

\newcommand{\bFormula}[1]{ \begin{equation} \label{#1} }
	\newcommand{\eF}{ \end{equation} }

\newcommand{\vr}{\varrho}

\newcommand{\vu}{\vc{u}}
\newcommand{\vc}[1]{{\bf #1}}

\newcommand{\Div}{{\rm div}_x}
\newcommand{\Grad}{\nabla_x}

\newcommand{\dx}{{\rm d} {x}}
\newcommand{\dy}{{\rm d} {y}}
\newcommand{\dt}{{\rm d} t }

\newcommand{\vU}{\vc{U}}

\newcommand{\intO}[1]{\int_{\mathbb{T}^d} #1 \dx}

\font\F=msbm10 scaled 1000
\newcommand{\R}{\mbox{\F R}}

\definecolor{Cgrey}{rgb}{0.85,0.85,0.85}
\definecolor{Cblue}{rgb}{0.50,0.85,0.85}
\definecolor{Cred}{rgb}{1,0,0}
\definecolor{fancy}{rgb}{0.10,0.85,0.10}

\newcommand\Cbox[2]{%
	\newbox\contentbox%
	\newbox\bkgdbox%
	\setbox\contentbox\hbox to \hsize{%
		\vtop{
			\kern\columnsep
			\hbox to \hsize{%
				\kern\columnsep%
				\advance\hsize by -2\columnsep%
				\setlength{\textwidth}{\hsize}%
				\vbox{
					\parskip=\baselineskip
					\parindent=0bp
					#2
				}%
				\kern\columnsep%
			}%
			\kern\columnsep%
		}%
	}%
	\setbox\bkgdbox\vbox{
		\color{#1}
		\hrule width  \wd\contentbox %
		height \ht\contentbox %
		depth  \dp\contentbox
		\color{black}
	}%
	\wd\bkgdbox=0bp%
	\vbox{\hbox to \hsize{\box\bkgdbox\box\contentbox}}%
	\vskip\baselineskip%
}

\date{}

\begin{document}
	\title{ Dissipative measure-valued solutions and weak-strong uniqueness \\for the Euler alignment system}
\author[1]{A. Chaudhary\thanks{Email: \href{mailto:chaudhary@na.uni-tuebingen.de}{chaudhary@na.uni-tuebingen.de}}}
	\author[2]{U. Koley\thanks{Email: \href{mailto: ujjwal@tifrbng.res.in}{ujjwal@math.tifrbng.res.in}}}

	\author[3]{E. Wiedemann\thanks{Email: \href{mailto:emil.wiedemann@fau.de}{emil.wiedemann@fau.de}}}
	
	\affil[1]{Universität Tübingen, Auf der Morgenstelle 10, 72076 Tübingen, Germany}
	
	\affil[2]{TIFR-CAM, P.O. Box 6503, GKVK Post Office, 560065, Bangalore, India }
	
	\affil[3]{Friedrich-Alexander-Universit\"at Erlangen-Nürnberg, Cauerstr.~11, 91058 Erlangen, Germany}
	\maketitle
	\begin{abstract}
	We introduce the concept of a dissipative measure-valued solution to the Euler alignment system. This approach incorporates a modified total energy balance, utilizing a binary tensor Young measure. The central finding is a weak (measure-valued)--strong uniqueness principle: if both a dissipative measure-valued solution and a classical smooth solution originate from the same initial data, they will be identical as long as the classical solution exists.
	\end{abstract}
	


	\section{Introduction}
	\label{I}
	The study of hydrodynamic models originating from the field of self-organization has garnered significant attention in recent years, leading to new challenges in partial differential equations that pose intriguing questions for the mathematical community. In this paper, we are interested in the following Euler alignment system
	\begin{equation}\label{NS}
	\begin{cases}
				\partial_t \vr + \Div (\vr \vu )=0\qquad&\text{in}\,\mathbb{T}^d\times(0,T),\\
				\partial_t (\vr \vu) + \Div (\vr \vu \otimes \vu) + \Grad p(\vr)= \mathcal{L}^\lambda[\vr, \vu]\qquad&\text{in}\,\mathbb{T}^d\times(0,T),\\
				(\vr(0), \vu (0))=(\vr_0, \vu_0)\qquad&\text{in}\,\mathbb{T}^d,
			\end{cases}
	\end{equation}
where $\lambda \in (0,1)$ is a fixed number, $\vr$ is the density, $\vu$ is the velocity, $p(\vr)=\vr^\gamma$ is the given pressure function with $\gamma\ge2$. Moreover, $\mathbb{T}^d$ denotes the $d$-dimensional periodic torus $[0, 2\pi]^d$, with endpoints identified for $d =1,2,3$. Finally, the commutator $\mathcal{L}^\lambda[\vr,\vu]$ is defined as
\begin{align}\label{commutator}
		\mathcal{L}^\lambda[\vr,\vu](t,x):=\vr(t,x) \mathfrak{L}_{\lambda}[\vr(t,\cdot) \vu(t,\cdot)](x) - \vr(t,x) \vu(t,x) \fr[\vr(t,\cdot)](x),
\end{align}
where the fractional Laplacian $\fr[\cdot]$ of order $\lambda\in (0,1)$ defined on the torus $\mathbb{T}^d$ as 
\begin{align*}
		\fr[\varphi](x)= [(-\Delta)^\lambda{\mathbf{ \varphi}}](x)
		:=- \text{P.V}\int_{\mathbb{T}^d}({\mathbf{ \varphi}}(x)-{\mathbf{ \varphi}}(y))\kappa_\lambda(x-y)dy,
	\end{align*}
for a sufficiently regular function $\varphi$. Here the singular kernel of the fractional Laplacian given by 
$$\kappa_\lambda(x):=\sum_{n\in\mathbb{Z}^d}\frac{1}{|x-2n\pi|^{d+2\lambda}}, \qquad 0<\lambda<1.
$$ 
The Euler alignment system \eqref{NS} is a model that effectively captures the collective behavior observed in animal swarms. In the past decade, this system has garnered increasing attention in the literature, particularly concerning its global well-posedness and the long-term flocking behavior, see \cite{Tadmor&Tan}.
	 
	 The commutator \eqref{commutator} with the singular kernel contributes to the system's dissipative nature. Specifically, when enforcing $\vr \equiv 1$, the alignment term transforms into a fractional Laplacian, which induces a regularizing effect on the solution. The corresponding system
	 \begin{align}\label{fractional Euler equation}
	 	\partial_t \vu + \vu\cdot\nabla \vu= (-\Delta)^{\lambda}\vu
	 \end{align}
 is recognized as the fractal Burgers equation. This equation has been comprehensively analyzed in \cite{Kiselev&Nazarov&Shterenberg}, where global well-posedness is demonstrated for $\lambda\in [1/2,1)$. In contrast, for $\lambda\in (0,1/2)$, solutions are prone to forming shocks. The regularization effect of $\eqref{commutator}$ was investigated in $\cite{Shvydkoy&Tadmor2017, Shvydkoy&Tadmor2018, Do&Kiselev&Ryzhik&Tan2018}$ for a one-dimensional periodic domain  $\mathbb{T}$, assuming the absence of pressure. In \cite{Do&Kiselev&Ryzhik&Tan2018}, a notable finding is the global well-posedness of the system for any 
 $\lambda\in(0,1)$. Specifically, when $\lambda\in (0,1/2)$, unlike in $\eqref{fractional Euler equation}$, the regularization effect of $\eqref{commutator}$ is sufficiently strong to guarantee global regularity for all smooth initial data.
 
 The theory of the pressureless Euler alignment system becomes increasingly complex in higher dimensions. Global well-posedness has predominantly been established for small initial data that are perturbed around \( \rho \equiv 1 \). Notable recent results include Shvydkoy's work \cite{Shvydkoy2019}, which addresses the periodic domain \( \mathbb{T}^d \) with \( \lambda \in (0,1) \), and the study by Danchin et al.\ \cite{Denchin&Mucha&Peszek&Bartosz2019}, which focuses on the entire space \( \mathbb{R}^d \) with \( \lambda \in (1/2,1) \).  Global regularity for general initial data is proven for uni-directional flows \cite{Lear2023, Lear&Shvydkoy2021, Li&Miao&Tan&Xue2023} with \( \lambda \in [1/2,1) \). 
 
  The global well-posedness theory for the Euler-alignment system \eqref{NS} with pressure is much less understood compared to the pressureless system. For the one-dimensional torus \( \mathbb{T} \), Constantin et al. \cite{Constantin&Drivas&Shvydkoy2020} established a global well-posedness theory for sufficiently large \( \lambda \in (5/6,1) \). Their result necessitates the inclusion of an additional strong local dissipation term of the form \( (\mu(\rho)u_x)_x \) in the system. In higher dimensions, existing global regularity results typically depend on smallness conditions applied to the initial data. For \( \lambda \in (0,1) \), Chen et al.\ \cite{Chen&Tan&Tong2021} demonstrated global well-posedness for smooth initial data under smallness assumptions in the spatial domain \( \mathbb{T}^d \). 
  
   For the Euler-alignment system \eqref{NS} in \( \mathbb{R}^d \), the authors in \cite{Bai&Miao&Tan&Xue2024} investigated global well-posedness for small initial data within a suitable critical Besov space, specifically when \( \lambda \in (1/2,1) \). Recently, in \cite{Bai&Tan&Liutang2024} the authors studied the Cauchy problem of the compressible Euler system with strongly singular velocity alignment. They established a global well-posedness theory for the system with small smooth initial data for $\lambda\in (0,1/2]$.  Global regularity for general large initial data remains a challenging open problem.
   
B\v{r}ezina and M\'{a}cha analyzed the equations governing a viscous approximation of the generalized compressible Euler system in \cite{Brezina&Macha}. They demonstrated that, as viscosity approaches zero, the dissipative measure-valued solutions converge to a strong solution of the Euler system, provided such a strong solution exists. In recent work \cite{2024}, Carrillo et al. studied several pressureless variants of the compressible Euler equation driven by nonlocal repulsion-attraction and alignment forces with Poisson interaction. Under an energy admissibility criterion, they established the existence of global measure-valued solutions.
   
 \subsection{ Scope and outline of the paper}
 
 The concept of measure-valued solutions to partial differential equations was introduced by DiPerna in the context of conservation laws (see \cite{DiPerna1985}). By employing Young measures, he successfully handled the passage to the artificial viscosity limit. In specific cases, such as the scalar setting, he also demonstrated a posteriori that these measure-valued solutions are atomic, implying they coincide with distributional solutions. For general systems of conservation laws, obtaining distributional (or entropy) solutions is often impractical or impossible. Therefore, measure-valued solutions or related frameworks become essential. In fluid dynamics, the existence of measure-valued solutions has been established for various models (\cite{DiPerna&Majda1987, Gwiazda2005, Neustupa1993, Emil_NSE}), highlighting their significance in addressing complex systems.

We aim to investigate the Euler alignment system \eqref{NS} with pressure for generic initial data. We develop a suitable mathematical framework for {\em {\em dissipative measure-valued} } solutions to the Euler alignment system \eqref{NS}, with the goal of establishing a weak (measure-valued)–strong uniqueness principle. Notably, our solution framework relies solely on natural energy bounds and a uniform lower bound on density related to approximate solutions. 
    
In this work, we focus on measure-valued solutions because they can capture a broad class of approximate problems, including systems with higher-order viscosities and specific numerical schemes. These approaches naturally produce measure-valued solutions, even in cases where proving convergence to a weak solution is difficult or remains unresolved. Motivated by this, we introduce the notion of a (dissipative) measure-valued solution for equation \eqref{NS}.
  
The key novelty of our work is the successful identification, in the limit, of the non-local term \eqref{commutator}, which is crucial for the framework of {\em {\em dissipative measure-valued}} solutions (see Definition \ref{DD1}) and proving the weak (measure-valued)–strong uniqueness principle (see Theorem \ref{second result}). In the case of a non-local term \eqref{commutator}, the standard approach of identifying non-linear terms using a Young measure is  not directly effective. To address this, we used the concept of a binary tensor Young measure, see Definition \ref{def1}, which ensures the accurate identification of the non-local term in the limit. It turns out that establishing a compatibility condition (see \eqref{compatibility second}) and the core proprieties of a binary tensor Young measure, which are inherently satisfied by any measure derived from a well-behaved sequence of approximate solutions to \eqref{NS} (see Theorem \ref{first result}), is enough to ensure weak-strong uniqueness. As a result, we bypass the cumbersome notation used in the framework of Alibert and Bouchitté \cite{Alibert&Bouchitte1997} and provide a more comprehensive definition of {\em {\em dissipative measure-valued}} solutions that still adheres to the weak-strong uniqueness principle.
    
This paper's central objective is to establish weak-strong uniqueness for dissipative measure-valued solutions to \eqref{NS} as defined in Definition \ref{DD1}. Weak-strong uniqueness guarantees that classical solutions are stable within the more general framework of dissipative measure-valued solutions. Notably, we present the first proof of weak-strong uniqueness for measure-valued solutions of the Euler alignment model \eqref{NS}.
     
We also identify a wide range of problems that generate dissipative measure-valued solutions, including density-dependent viscosity equations, which are outside the scope of current weak solution theory. In addition, we observe that certain numerical schemes are consistent with this approach, aligning with the viewpoint of Fjordholm et al. \cite{Fjordholm&Kappeli&Mishra&Tadmor2014}. They contend that for hyperbolic systems of conservation laws, dissipative measure-valued solutions provide a more suitable framework than weak entropy solutions. This is because dissipative measure-valued solutions can be derived as limits of typical numerical approximations, a feat not possible with weak entropy solutions.  
        
The rest of this paper is organized as follows: Section~\ref{Section 2} presents the key mathematical and technical frameworks, followed by a discussion of solution concepts and the main results. Section~\ref{Section 3} demonstrates the convergence of approximate solutions and proves the existence of {\em dissipative measure-valued}  solutions to \eqref{NS}. Finally, in Section~\ref{Section 4}, we establish a weak (measure-valued)-strong uniqueness principle for \eqref{NS}, based on a suitable relative energy inequality for the Euler alignment system.\\

{\bf Acknowledgement.} E.~W.~is supported by DFG grant no.~525716336 within the Priority Programme SPP 2410.

\section{Mathematical framework and main results}
\label{Section 2}
We start with a brief review of the key mathematical tools essential for the subsequent analysis, followed by a presentation of the paper's main results. Throughout, the symbol C denotes generic constants that may vary from one step to another in the proofs. While it is theoretically possible to track these constants explicitly, this would be overly cumbersome, so we omit such details for clarity.

	\subsection{Weak solution to  the regularized system}
	We will approximate a {\em {\em dissipative measure-valued}} solution to \eqref{NS} as a limit of viscous dissipative weak solutions to the following regularized system:
	\begin{equation}\label{NS1}
		\begin{aligned}
			\begin{cases}
				\partial_t \vr_\varepsilon + \Div (\vr_\varepsilon \vu_\varepsilon)=0\qquad&\text{in}\,\mathbb{T}^d\times[0,T],\\
				\partial_t (\vr_\varepsilon \vu_\varepsilon) + \Div (\vr_\varepsilon \vu_\varepsilon \otimes \vu_\varepsilon) + \Grad p(\vr_\varepsilon)= \mathcal{L}^\lambda[\vr_\varepsilon,\vu_\varepsilon]-\varepsilon \Delta^{2m} \vu_\varepsilon\qquad&\text{in}\,\mathbb{T}^d\times[0,T],\\
				(\vr_\varepsilon(0), \vu_\varepsilon( 0))=(\vr_{0,\varepsilon}, \vu_{0,\varepsilon})\qquad&\text{in}\,\mathbb{T}^d,
			\end{cases}
		\end{aligned}
	\end{equation}
	where $m\in\mathbb{N}$, $(\vr_{0,\varepsilon}, \vu_{0,\varepsilon})$ is given by $\vr_{0,\varepsilon}=\vr_0 *g_\varepsilon $, $\vu_{0,\varepsilon}=\vu* {g}_{\varepsilon}$, and $g_\varepsilon (\cdot)=\frac{1}{\varepsilon^d}g(\frac{\cdot}{\varepsilon})$ with a standard mollifier $g$ on the torus $\mathbb{T}^d$.
	
For sufficiently large $m\in\mathbb{N}$ and fixed $\varepsilon>0$, modulo cosmetic changes, the existence of a unique smooth solution to the regularized system \eqref{NS1} can be established; see, for example, \cite{Kroner, BrFeMa2020} for a detailed discussion. Next, we introduce the weak formulation associated with the regularized system \eqref{NS1}. This is essential because, in the final steps of proving the existence of a dissipative measure-valued solution, we will need to pass to the limit in this weak formulation.
	\begin{Definition}\label{dissipative weak solution}
	Let $(\vr_0, \vu_0)$ be an initial data such that $\intO{ \left( \frac{1}{2} \vr_0 |\vc{u}_0|^2 + P(\vr_0) \right)}\,\textless\,\infty$. We call $(\vr_\varepsilon, \vu_\varepsilon)\in L^\infty([0,T];L^{\gamma}(\mathbb{T}^d))\times L^2([0,T];H^{2m}(\mathbb{T}^d))$ a dissipative weak solution to \eqref{NS1} if $(\vr_\varepsilon, \vu_\varepsilon)$ satisfies the following conditions:
	\begin{itemize}
		\item {\bf Equation of continuity.}
		\begin{equation} \label{viscouse continuity equation}
			\begin{split}
				&\intO{ \vr_\varepsilon(\tau, x) \psi (\tau, x) } -  \intO{ \vr_{0,\varepsilon}(x) \psi (0,x) } \\
				& \qquad = \int_0^\tau \int_{\mathbb{T}^d}{ \Big[ \vr_\varepsilon(t, x) \partial_t \psi(t,x) + \vr_\varepsilon(t,x) \vu_\varepsilon(t, x) \cdot \Grad \psi(t,x) \Big] }\dx\dt,
			\end{split}
		\end{equation}
		holds for every $\psi \in C^1([0,T] \times \mathbb{T}^d)$ and $\tau\in[0,T]$.
		\item {\bf Momentum equation.}
		\begin{align}
		&\intO{ \vr_\varepsilon(\tau,x) \vu_\varepsilon(\tau,x) \cdot {\mathbf{ \varphi}} (\tau, x) }  -  \intO{ \vr_{0,\varepsilon}(x)\vu_{0,\varepsilon}(x) \cdot {\mathbf{ \varphi}} (0, x) } + \varepsilon\int_0^\tau \int_{\mathbb{T}^d} \Delta^{m}\vu_{\varepsilon}(t,x)\cdot\Delta^m\varphi(t,x){\rm d}x{\rm d}t \nonumber \\
				&= \int_0^\tau \int_{\mathbb{T}^d}{ \Big[ \vr_\varepsilon(t,x) \vu_\varepsilon(t,x) \cdot \partial_t {\mathbf{ \varphi}(t,x)} + \vr_\varepsilon(\vc{u}_\varepsilon (t,x)\otimes \vc{u}_\varepsilon(t,x)) :\Grad {\mathbf{ \varphi}(t,x)} + p(\vr_\varepsilon(t,x)) \Div {\mathbf{ \varphi}(t,x)} \Big]}\dx\dt \nonumber \\
				& \quad - \frac 12\int_0^\tau \int_{\mathbb{T}^{2d} }  \vr_\varepsilon(t,x) \vr_\varepsilon(t,y){(\vc{u}_\varepsilon(t,x) - \vc{u}_\varepsilon(t,y)) ({\mathbf{ \varphi}}(t,x)-{\mathbf{ \varphi}}(t,y))}\kappa_\lambda(x-y)\dx\dy\dt, \label{viscous momentum equation}
			\end{align}
		holds for every ${\mathbf{ \varphi}} \in C^1([0,T]; H^{2m} (\mathbb{T}^d; \R^d))$ and $\tau\in[0,T]$.
		\item{\bf Energy inequality.}
		\begin{equation}\label{viscouse energy inequality}
			\begin{split}
				&\int_{\mathbb{T}^d}\left( \frac{1}{2} \vr_\varepsilon |\vc{u}_\varepsilon(\tau,x)|^2 +P(\vr_\varepsilon(\tau,x))\right){\rm d}x +\varepsilon\int_0^\tau\int_{\mathbb{T}^d}|\Delta^{m}\vu_\varepsilon(t,x)|^2{\rm d}x {\rm d}t
				\\&\qquad\qquad\qquad
				+\frac{1}{2}\int_0^\tau \int_{\mathbb{T}^{2d}} \vr_\varepsilon(t, x) \vr_\varepsilon(t, y){|\vc{u}_\varepsilon(t, x) - \vc{u}_\varepsilon(t,y)|^2}\kappa_\lambda(x-y) \dx\dy\dt \\
				&\qquad \leq \intO{ \left( \frac{1}{2} \vr_{0,\varepsilon}(x) |\vc{u}_{0,\varepsilon}(x)|^2 +\ P(\vr_{0,\varepsilon}(x)) \right)},
			\end{split}
		\end{equation}
		for a.a. $\tau \in (0,T)$, where the pressure potential $P$ is given by $P(\vr)=\frac{\vr^\gamma}{\gamma-1}$.
	\end{itemize}
\end{Definition}

\subsection{Binary tensor Young measure}{Given the absence of well-posedness results for weak solutions to \eqref{NS}, there is significant potential for developing alternative solution frameworks for equation \eqref{NS}. One such framework involves measure-valued solutions, where the solutions sought are not functions but rather space-time parameterized probability measures on the state space. Motivated by the concept of statistical solutions \cite{Fjordholm&Kappeli&Mishra&Tadmor2014, Fjordholm&Lanthaler&Mishra2017}, we introduce the notion of a binary tensor Young measure. Given a topological space $X$, we denote by $\mathcal{P}(X)$ the space of all probability measures on the Borel sigma-algebra of $X$. We denote by $\mathcal{M}(X)$ (resp.\ $\mathcal{M}^+(X)$) the space of all measures (resp.\ positive measures) on the Borel sigma-algebra of $X$. For convenience we set $U=  \mathbb{R}^d\times[0,\infty)$.
\begin{Definition}[Binary tensor Young measure] \label{def1} 
	Let \(\nu: [0,T] \times \mathbb{T}^d \to \mathcal{P}(U)\) be a weak*-measurable mapping, meaning that for every \(f \in C_0(U)\), the function \((t,x) \mapsto \langle \nu_{t,x}, f \rangle\), mapping \([0,T] \times \mathbb{T}^d\) into \(\mathbb{R}\), is Borel measurable. In other words, \(\nu\) is a Young measure. The binary tensor Young measure, denoted by \(\nu^2_{t,x,y}\), is defined as the product measure \(\nu_{t,x} \otimes \nu_{t,y}\).
\end{Definition}

It is easy to see that, if $(u_n)_{n\in\mathbb N} \in L^1(\mathbb{T}^{d};\R^k)$ generates the Young measure $\nu$, then the sequence given by $(x,y)\mapsto (u_n(x),u_n(y))\in (\R^k)^2$ generates the tensorized Young measure $\nu^2_{x,y}(dz_1dz_2)=\nu_x(dz_1)\nu_y(dz_2)$. We now collect some useful properties of binary tensor Young measure.

\begin{Remark}[Properties of binary tensor Young measure] 
\label{remark1} 
A binary tensor Young measure $\nu_{t,x,y}^2=\nu_{t,x}\otimes\nu_{t,y}$ satisfies the following properties:
	\begin{itemize}
				\item[1.] \textbf{Symmetry:} If $f\in C_0(U^2)$, then
				\begin{align*}
					\left<\nu_{t,x,y}^2; f(s, \vc{v}, s', \vc{v}')\right> = \left<\nu_{t,y,x}^2; f(s', \vc{v}', s, \vc{v})\right>
				\end{align*}
				for a.e. $(t,x,y)\in [0,T]\times\mathbb{T}^d\times\mathbb{T}^d$.
				\item[2.] \textbf{Consistency:} If $f\in C_0(U^2)$ is of the form $f(s,\vc{v},s',\vc{v}')= g(s,\vc{v})$, for some $g\in C_0(U)$, then 
				\begin{align*}
					\left<\nu_{t,x,y}^2; f\right>=\left<\nu_{t,x}; g\right>
				\end{align*}
				for a.e. $(t,x,y)\in [0,T]\times\mathbb{T}^d\times\mathbb{T}^d$.
			\end{itemize}
			In view of the above properties it is clear that $\nu_{t,x,y}^2$ contains all informations about $\nu_{t,x}$, but not  vice-versa.
	\end{Remark}}
	
%
%
%

\subsection{{{Dissipative measure-valued}} solution}
Now we are in a position to introduce the concept of a {\em dissipative measure-valued solution} of the Euler alignment system \eqref{NS}. 
	\begin{Definition} \label{DD1}
		We say that a pair $[ \nu, \mathcal{D} ]$
		is a dissipative measure-valued solution of the Euler alignment system \eqref{NS}, in $(0,T) \times \mathbb{T}^d$ with the given initial data $\nu_0$, if the following conditions hold:
		\begin{itemize}
			\item {\bf Young measure.}
		  \(\nu: [0,T] \times \mathbb{T}^d \to \mathcal{P}(U)\) is a weak*-measurable mapping, meaning that for every \(f \in C_0(U)\), the function \((t,x) \mapsto \langle \nu_{t,x}, f \rangle\), mapping \([0,T] \times \mathbb{T}^d\) into \(\mathbb{R}\), is Borel measurable. In other words, \(\nu\) is a Young measure.
		 
		 \item {\bf{$L^\infty$-boundedness}.} $\nu_{t,x}$ satisfies the following bounds
		 \begin{align}\label{boundedness}
		 	\sup_{t\in[0,T]}\int_{\mathbb{T}^d}\left<\nu_{t,x};|s|^\gamma+ s|\vc{v}|^{2}\right>dx\,<\,\infty.
		 \end{align}
	 
			\item {\bf Dissipation defect.} The dissipation defect satisfies
		\[
		\mathcal{D} \in L^\infty(0,T), \ \mathcal{D} \geq 0.
		\]
	
			\item {\bf Equation of continuity.}
			There exists a measure $r^C\in L^1([0,T];\mathcal{M}(\mathbb{T}^d))$ such that for a.a.\ $\tau\in (0,T)$ and for every $\psi \in C^1([0,T] \times \mathbb{T}^d)$,
			\begin{equation}
				\left| \langle r^C (\tau) ; \Grad\psi \rangle \right| \leq \chi(\tau)  \mathcal{D} (\tau) \| \psi \|_{C^1({\mathbb{T}^d})}
			\end{equation}
		and
			{\begin{equation} \label{dmvB1}
					\begin{split}
						\int_{\mathbb{T}^d}{ \langle \nu_{\tau,x}; s \rangle \psi (\tau, \cdot) } {\rm d}x &-  \int_{\mathbb{T}^d}{ \langle (\nu_0)_x; s \rangle \psi (0, \cdot) }{\rm d}x \\
						&= \int_0^\tau \int_{\mathbb{T}^d}{ \Big[ \langle \nu_{t,x}; s \rangle \partial_t \psi (t,x) + \langle \nu_{t,x}; s\vc{v} \rangle \cdot \Grad \psi(t,x) \Big] } {\rm d}x \dt
						+ \int_0^\tau \langle r^C; \Grad \psi \rangle \ \dt.
					\end{split}
			\end{equation}}
			
			\color{black}
			
			\item {\bf Momentum equation.}
			There exists a measure $r^M\in L^1([0,T];\mathcal{M}(\mathbb{T}^d))$ such that for a.a.\ $\tau\in(0,T)$ and for every ${\mathbf{ \varphi}} \in C^1([0,T] \times \mathbb{T}^d; \mathbb{R}^d)$, 
			\begin{equation}\label{concentrate defect 2}
				\left| \langle r^M (\tau) ; \Grad\varphi \rangle \right| \leq \xi(\tau)  \mathcal{D} (\tau) \| \varphi \|_{C^1({\mathbb{T}^d})}
			\end{equation}
			and
			\begin{equation} \label{dmvB2}
				\begin{split}
					&\int_{{\mathbb{T}^d}}{ \langle \nu_{\tau,x};  s\vc{v} \rangle \cdot {\mathbf{ \varphi}} (\tau, \cdot) } dx -  \int_{\mathbb{T}^d}{ \langle (\nu_0)_x;s\vc{v} \rangle \cdot {\mathbf{ \varphi}} (0, \cdot) }dx \\
					&= \int_0^\tau \int_{\mathbb{T}^d}{ \Big[ \langle \nu_{t,x} ; s\vc{v} \rangle \cdot \partial_t {\mathbf{\varphi}(t,x)} + \langle \nu_{t,x}; {s\vc{v} \otimes \vc{v}}\rangle : \Grad {\mathbf{ \varphi}(t,x)} +
						\langle \nu_{t,x}; p(s) \rangle \Div {\mathbf{ \varphi}(t,x)} \Big] }dx\ \dt\\
					& -\frac{1}{2}{\int_0^\tau \int_{\mathbb{T}^d\times\mathbb{T}^d} \left< \nu_{t,x,y}^2 ; ss'{(\vc{v}-\vc{v}') ({\mathbf{ \varphi}}(t,x)-{\mathbf{ \varphi}}(t,y))}{\kappa_\lambda(x-y)}\right> \dx\dy\dt }+ \int_0^\tau \left< {r}^M ; \Grad {\mathbf{ \varphi}} \right>\dt.
				\end{split}
			\end{equation}
			
			\color{black}
			
			\item{\bf Energy inequality.} The following energy inequality holds:
			\begin{equation} \label{dmvEI}
				\begin{split}
					\int_{\mathbb{T}^d}{ \left< \nu_{\tau,x};  \left( \frac{1}{2} {s|\vc{v}|^2} + P(s) \right) \right> }\dx\dy
					&+\frac{1}{2}{\int_0^\tau \int_{\mathbb{T}^d\times \mathbb{T}^d} \left<\nu_{t,x,y}^2; {ss'|\vc{v}-\vc{v}'|^2} \right> \kappa_\lambda(x-y)\dx\dy\dt}\\+\mathcal{D}(\tau)
					&\leq \int_{\mathbb{T}^d}{ \left< (\nu_{0})_x; \left( \frac{1}{2}{s|\vc{v}|^2}+ P(s) \right) \right>}\dx
				\end{split}
			\end{equation}
			for a.a. $\tau \in (0,T)$.
			\item {\bf Compatibility conditions.} If $\mathbf{U}\in C([0,T];C^1(\mathbb{T}^d))$, then the following condition holds:
			\begin{align}\label{compatibility second}
				&\int_0^\tau\int_{\mathbb{T}^d\times\mathbb{T}^d} \left<\nu_{t,x}; |\big(\vc{v}-\bar{\vc{v}}(t)\big)-\big(\mathbf{U}(t,x)-\bar{\vU}(t)\big)|^2\right>\dx\dy\dt\notag\\&\le C_\lambda\int_0^\tau\int_{{\mathbb{T}^d}\times\mathbb{T}^d}\left<\nu_{t,x,y}^2; ss'|(\vc{v}-\vc{v}')-(\mathbf{U}(t,x)-\mathbf{U}(t,y))|^2 \kappa_\lambda(x-y)\right>\dx\dy\dt+C_\lambda\mathcal{D}(\tau)
			\end{align}
		 for a.a.\ $\tau\in [0,T]$, for some positive constant $C_\lambda$, where $\bar{\vc{v}}(t)=\int_{\mathbb{T}^d}\left<\nu_{t,x},\vc{v}\right>{\rm d}x$ and $\bar{\vU}(t)=\int_{\mathbb{T}^d}\vU(t,x){\rm d}x$.
		\item {\bf Support of $\nu$.} There exists a constant $c_\vr>0$ such that
		\begin{align}\label{supp1}
			\mbox{supp}(\nu_{t,x,y}^2)\subset \{(s,\vc{v}, s',\vc{v}'): s,s'\ge c_\vr>0\,\&\, \vc{v},\vc{v}'\in\mathbb{R}^d\},
		\end{align}
		and 
		\begin{align}\label{supp2}
			\mbox{supp}(\nu_{t,x})\subset \{(s,\vc{v}): s\ge c_\vr>0\,\&\, \vc{v}\in\mathbb{R}^d\},
		\end{align}
		for a.e.\ $(t,x,y)\in[0,T]\times\mathbb{T}^d\times\mathbb{T}^d.$
		\end{itemize}
	\end{Definition}
\begin{Remark}
	Note that \eqref{boundedness} and \eqref{supp2} give that
	\begin{align}\label{L^2 boundedness}
		\sup_{t\in[0,T]}\int_{\mathbb{T}^d}\left<\nu_{t,x};|\vc{v}|^2\right>{\rm d}x < \infty.
	\end{align}
\end{Remark}
\subsection{Main results}
We now present the main results of this paper. First, concerning the convergence of the regularized system \eqref{NS1} to a dissipative weak solution to \eqref{NS}, we have the following theorem.
\begin{Theorem}[\textbf{Convergence to a {\em {\em dissipative measure-valued}} solution}]\label{first result} Let $(\vr_0, \vu_0)$ be an initial data such that $\intO{ \left( \frac{1}{2} \vr_0 |\vc{u}_0|^2 + P(\vr_0) \right)}\,\textless\,\infty$. Let $\{(\vr_\varepsilon,\vu_\varepsilon)\}_{\varepsilon>0}$ be a family of viscous dissipative weak solutions to \eqref{NS1} in the sense of Definition \ref{dissipative weak solution} and assume that there exists a $\underline{\vr}>0$ such that $\underline{\vr}\le\,\vr_\varepsilon$ for all $\varepsilon>0$. Then the family $\{(\vr_{\varepsilon},\vu_{\varepsilon})\}_{\varepsilon>0}$ generates a {\em {\em dissipative measure-valued}} solution $[ \nu, \mathcal{D} ]$ to \eqref{NS} in the sense of Definition \ref{DD1}.
\end{Theorem}
We then establish the following weak (measure-valued)-strong uniqueness principle.
\begin{Theorem}[\textbf{Weak (measure-valued)-strong uniqueness principle}] \label{second result}
	Let $[ \nu, \mathcal{D} ]$ be a {dissipative measure-valued} 
	solution to the Euler alignment system \eqref{NS}, as described in Definition \ref{DD1}.
	Let $(r, \vc{U})$ be a smooth solution of the Euler alignment system, satisfying the following 
	\begin{align}\label{assumption}
		r, \ \Grad r, \ \vc{U},\ \Grad \vc{U} \in C([0,T] \times {\mathbb{T}^d}),\
		\partial_t \vc{U} \in L^\infty(0,T; C({\mathbb{T}^d};\R^d)),\ r > 0.
	\end{align}
	Then there is a constant $\Sigma = \Sigma(T)$, depending on functions $r$ and $\vc{U}$, such that for a.a.\ $\tau \in (0,T)$ the following inequality holds:
	\[
	\begin{split}
		& \intO{ \left< \nu_{\tau,x};  \frac{1}{2} s |\vc{v} - \vc{U}(\tau,x)|^2 + P(s) - P'(r(\tau,x)) (s - r(\tau,x)) - P(r(\tau,x)) \right> }
		\\ &+\int_0^\tau \int_{\mathbb{T}^d \times \mathbb{T}^d} \left< \nu_{\tau,x,y}^2; s s' {|(\vc{v} - \vc{v}')- (\vc{U}(\tau,x)-\vc{U}(\tau,y))|^2}\right>\kappa_\lambda(x-y) \dx\dy\dt + \mathcal{D}(\tau) \\
		&\leq \Sigma(T) \intO{ \left< (\nu_{0})_x;  \frac{1}{2} s |\vc{v} - \vc{U}(0, \cdot) |^2 + P(s) - P'(r(0,\cdot)) (s - r(0,\cdot)) - P(r(0,\cdot)) \right> }
	\end{split}
	\]
This means, if	
\[
	(\nu_0)_x = \delta_{[ r(0,x), \vc{U}(0,x) ]}, \ \mbox{for a.a.} \ x \in \mathbb{T}^d,
	\]
	then $\mathcal{D} = 0$ and
	\[
	\nu_{\tau,x} = \delta_{[ r(\tau,x), \vc{U}(\tau,x) ]}, \ \mbox{for a.a.}\ \tau \in (0,T),\ x \in \mathbb{T}^d.
	\]
	
\end{Theorem}
\section{Proof of Theorem \ref{first result}}
\label{Section 3}
Let $\{(\vr_\varepsilon, \vu_\varepsilon)\}_{\varepsilon>0}$ be a family of viscous dissipative weak solutions to \eqref{NS1} with initial data $(\vr_{0,\epsilon}, \vu_{0,\epsilon})$ in the sense of Definition \ref{dissipative weak solution} such that the family  $\{(\vr_\varepsilon, \vu_\varepsilon)\}_{\varepsilon>0}$  satisfies the energy inequality 
		\begin{align}
			&\intO{\left( \frac{1}{2} \vr_{\varepsilon}(\tau, x) |\vc{u}_{\varepsilon}(\tau, x)|^2 + P(\vr_{\varepsilon}(\tau,x))\right)}
			+\frac 12 \int_0^\tau \int_{\mathbb{T}^d \times \mathbb{T}^d} \vr_{\varepsilon}(t,x) \vr_{\varepsilon}(t,y) {|\vc{u_{\varepsilon}}(t,x) - \vc{u}_{\varepsilon}(t,y)|^2}\kappa_\lambda(x-y)\dx\dy\dt  \nonumber \\
			&\qquad\qquad\qquad\qquad+\varepsilon\int_0^\tau\int_{\mathbb{T}^d}|\Delta^{m}\vu_\varepsilon(t,x)|^2 dx dt \leq \intO{ \left( \frac{1}{2} \vr_0(x) |\vc{u}_0(x)|^2 + P(\vr_0(x)) \right)}, \label{Energy estimate}
		\end{align}
and the following uniform lower bound on densities (by assumption)
\begin{align}\label{uniform lower bound}
	0<\underline{\vr}\le\,\vr_\varepsilon.
\end{align}
\subsection*{Identification of weak limits.}
 The energy estimate \eqref{Energy estimate} and uniform lower bound on densities \eqref{uniform lower bound} ensure that there exists a constant $C>0$ (independent of $\varepsilon$) such that
\begin{align}\label{energy bounds}
  \|\vr_\varepsilon\|_{L^\infty([0,T];L^\gamma(\mathbb{T}^d))}\le C,\quad
  \|\vr_\varepsilon\vu_\varepsilon\|_{L^\infty([0,T];L^{\frac{2\gamma}{\gamma+1}}(\mathbb{T}^d))}\,\le\,C,\quad \text{and}\quad\,\|\vu_{\varepsilon}\|_{L^\infty([0,T];L^2(\mathbb{T}^d))}\le\,C.
\end{align}
From \eqref{NS1} and \eqref{Energy estimate}, we can also obtain
\begin{align} \label{energy bounds_01}
	\|\vr_\varepsilon\|_{H^1([0,T];H^{-m}(\mathbb{T}^d))}\le\,C,\quad\text{and}\quad \|\vr_\varepsilon\vu_\varepsilon\|_{H^1([0,T];H^{-m}(\mathbb{T}^d))}\le\,C.
\end{align}
By using the standard weak compactness argument (Banach-Alaoglu theorem, Young measure compactness \cite{Blader} and compact embeddings), we can obtain that there exists a Young measure $\nu$ such that (up to a subsequence)
\begin{align*}
	&\nu_{t,x}^{1,\varepsilon}=\delta_{(\vr_\varepsilon(t,x), \vc{u}_{\varepsilon}(t,x))} \to \nu_{t,x}, \,\,\mbox{weak*}\,\text{in}\,L^\infty_{\rm weak}((0,T)\times \mathbb{T}^d;\mathcal{P}(U)),\\
&\vspace{2cm}\nu_{t,x,y}^{2,\varepsilon}=\delta_{(\vr_\varepsilon(t,x),\vr_{\varepsilon}(t,y), \vc{u}_\varepsilon(t,x),\vc{u}_{\varepsilon}(t,y))}\to \nu_{t,x,y}^2, \,\,\mbox{weak*}\,\text{in}\,L^\infty_{\rm weak}((0,T)\times (\mathbb{T}^d)^2;\mathcal{P}(U^2)),
	\\&\vspace{2cm}\vr_\varepsilon(t,x) \to \langle \nu_{t, x} ; s\rangle \,\,\text{ in}\,\,C_{\text{weak}}([0,T];L^\gamma(\mathbb{T}^d)),\,
	\\&\vspace{2cm}\vr_{\varepsilon}(t,x)\vc{u}_\varepsilon(t,x)\to \langle\nu_{t,x}; s \vc{v}\rangle\,\,\text{in}\, C_{\text{weak}}([0,T];L^{\frac{2\gamma}{\gamma+1}}(\mathbb{T}^d)),
	\\&\vspace{2cm}\vc{u}_\varepsilon(t,x)\to \langle\nu_{t,x};  \vc{v}\rangle\,\,\text{weak-*  in}\,\,L^\infty([0,T];L^{2}(\mathbb{T}^d)),\\&
	{\varepsilon}\Delta^m\vu_{\varepsilon}\to 0\,\,\text{weakly in}\,\,L^2([0,T];L^2(\mathbb{T}^d)),
\end{align*}
{ where convergence in $C_{\text{weak}}([0,T];L^\gamma(\mathbb{T}^d))$ means that for all $\varphi\in (L^\gamma(\mathbb{T}^d))^*$, $$\left<\vr_\epsilon, \varphi\right> \to \left<\left<v_{\cdot,x}, s\right>,\varphi\right>\,\,\text{in}\,\,C([0,T]).$$}

Notice that, in view of the a-priori estimates given by \eqref{energy bounds}--\eqref{energy bounds_01}, the convergence of all the terms present in \eqref{viscouse continuity equation} and \eqref{viscous momentum equation}, except the non-local term in \eqref{viscous momentum equation}, are quite straightforward. Therefore, we will only focus on the identification of the non-local term.

\noindent
\subsection*{Identification of the non-local term in the limit.} 

To proceed, let us first fix a test function $\varphi\in C([0,T]; C^1(\mathbb{T}^d;\mathbb{R}^d))$. In what follows, we identify the limit of the non-local term in several steps as follows: 

\noindent
\textbf{Step 1.}
For convenience, for any $\mathtt{M}>1$, we define for all $x,y\in (0,2\pi)^d$
\begin{align*}
	\kappa_\lambda(x-y)&=: \kappa_\lambda^{0,{1}}(x-y) + \kappa_\lambda^{2, \mathtt{M}}(x-y) +\kappa_\lambda^{\mathtt{M},\infty}(x-y),
\end{align*}
where
\begin{align*}	\kappa_\lambda^{0,{1}}(x-y)&=\sum_{|n|\le 1 ,n\in\mathbb{Z}^d}\frac{1}{|x-y-2\pi n|^{d+2\lambda}},\\
	\kappa_\lambda^{2,\mathtt{M}}(x-y)&=\sum_{1<|n|\le \mathtt{M} ,n\in\mathbb{Z}^d}\frac{1}{|x-y-2\pi n|^{d+2\lambda}},\\
	\kappa_\lambda^{\mathtt{M},\infty}(x-y)&=\sum_{|n|\textgreater \mathtt{M},n\in\mathbb{Z}^d}\frac{1}{|x-y-2\pi n|^{d+2\lambda}}.
\end{align*}

\vspace{0.3cm}

\noindent
\textbf{Step 2.} Let $f:U^2\to \mathbb{R}^d$ be a continuous function such that 
$$\int_{[0,T]\times \mathbb{T}^{2d}}|f(\vr_\epsilon(t,x),\vc{u}_\varepsilon(t,x), \vr_{\varepsilon}(t,y), \vc{u}_\varepsilon(t,y))|^r \dx\dy\dt\le\,C,\,\text{for some}\, r>1.
$$
By using weak convergence to the Young measure, we obtain
\begin{align}\label{weak-convergence0}
f(\vr_\epsilon(t,x),\vc{u}_\varepsilon(t,x), \vr_{\varepsilon}(t,y), \vc{u}_\varepsilon(t,y))\to \left<\nu_{t,x,y}^2; f(s,\vc{v}, s', \vc{v}')\right> \text{weakly in}\,L^r ([0,T]\times\mathbb{T}^{2d}).
\end{align}
In particular, we can take $$f(\vr_\epsilon(t,x),\vc{u}_\varepsilon(t,x), \vr_{\varepsilon}(t,y), \vc{u}_\varepsilon(t,y))=\vr_\epsilon(t,x)\vr_\epsilon(t,y)(\vc{u}_{\varepsilon}(t,x)-\vc{u}_\varepsilon(t,y)).$$
The uniform estimate \eqref{Energy estimate} implies that there exists a positive constant $C$ such that for $r=\gamma\wedge\frac{2\gamma}{\gamma+1}$,
\begin{align*}
	\int_{[0,T]\times\mathbb{T}^{2d}}|f(\vr_\epsilon(t,x),\vc{u}_\varepsilon(t,x), \vr_{\varepsilon}(t,y), \vc{u}_\varepsilon(t,y))|^r\dx\dy\dt\le C.
\end{align*}
It shows that
\begin{align}\label{weak-convergence}
	\vr_\epsilon(t,x)\vr_\epsilon(t,y)(\vc{u}_{\varepsilon}(t,x)-\vc{u}_\varepsilon(t,y))\to \left< \nu_{t,x,y}^2; ss'(\vc{v}-\vc{v}')\right>\,\text{weakly in}\,L^r ([0,T]\times\mathbb{T}^{2d}).
\end{align}

\noindent
\textbf{Step 2.} In this step, for the given $\varphi\in C([0,T]; C^1(\mathbb{T}^d
\mathbb{R}^d))$, we identify the limit for the following term:
\begin{align}\label{today7}
	I_\varepsilon&=\frac{1}{2}\int_0^\tau \int_{\mathbb{T}^d \times \mathbb{T}^d}  \vr_\varepsilon(t,x) \vr_\varepsilon(t,y){(\vc{u}_\varepsilon(t,x) - \vc{u}_\varepsilon(t,y)) ({\mathbf{ \varphi}}(t,x)-{\mathbf{ \varphi}}(t,y))}\kappa_\lambda(x-y)\dx\dy\dt\notag\\&=\frac{1}{2} \int_0^\tau\int_{\mathbb{T}^d\times\mathbb{T}^d}\vr_\varepsilon(t,x)\vr_\varepsilon(t,y){(\vu_\varepsilon(t,x)-\vu_\varepsilon(t,y)){\big(\mathbf{ \varphi}}(t,x)-{\mathbf{ \varphi}}(t,y)\big)}\kappa_\lambda^{0,1}(x-y) dx dy dt\notag\\ &\qquad+\frac{1}{2} \int_0^\tau\int_{\mathbb{T}^d\times\mathbb{T}^d}\vr_\varepsilon(t,x)\vr_\varepsilon(t,y){(\vu_\varepsilon(t,x)-\vu_\varepsilon(t,y)){\big(\mathbf{ \varphi}}(t,x)-{\mathbf{ \varphi}}(t,y)\big)}\kappa_\lambda^{2,\mathtt{M}}(x-y) dx dy dt\notag \\&\qquad+\frac{1}{2} \int_0^\tau\int_{\mathbb{T}^d\times\mathbb{T}^d}\vr_\varepsilon(t,x)\vr_\varepsilon(t,y){(\vu_\varepsilon(t,x)-\vu_\varepsilon(t,y)){\big(\mathbf{ \varphi}}(t,x)-{\mathbf{ \varphi}}(t,y)\big)}\kappa_\lambda^{\mathtt{M},\infty}(x-y) dx dy dt\notag\\&
	:= I_{\varepsilon}^{0,1} + I_{\varepsilon}^{2,\mathtt{M}} + I_{\varepsilon} ^{\mathtt{M},\infty}.
\end{align}
Now in the following sub-steps, we deal with the above terms separately.

\noindent
\textbf{Step 2(a)} In this sub-step, we estimate the term $I_\varepsilon^{0,1}$ as follows:
\begin{align*}
	I_\varepsilon^{0,1}=&\frac{1}{2}\int_{0}^\tau \int_{\mathbb{T}^d\times \mathbb{T}^d}\vr_\varepsilon(t,x)\vr_\varepsilon(t,y)\frac{(\vu_\varepsilon(t,x)-\vu_\varepsilon(t,y))(\varphi(t,x)-\varphi(t,y))}{|x-y|^{d+2\lambda}}\kappa(x-y){\rm d}x {\rm d}y {\rm d}t\\&\qquad +\frac{1}{2}\sum_{|n|=1,n\in\mathbb{Z}^d} \int_{0}^\tau \int_{\mathbb{T}^d\times \mathbb{T}^d}\vr_\varepsilon(t,x)\vr_\varepsilon(t,y)\frac{(\vu_\varepsilon(t,x)-\vu_\varepsilon(t,y))(\varphi(t,x)-\varphi(t,y))}{|x-y-2\pi n|^{d+2\lambda}}{\rm d}x {\rm d}y {\rm d}t,
\end{align*}

By using the periodicity of functions on $\mathbb{T}^d$, we can prove that 
\begin{align*}
	&\sum_{|n|=1,n\in\mathbb{Z}^d} \int_{0}^\tau \int_{\mathbb{T}^d\times \mathbb{T}^d}\vr_\varepsilon(t,x)\vr_\varepsilon(t,y)\frac{(\vu_\varepsilon(t,x)-\vu_\varepsilon(t,y))(\varphi(t,x)-\varphi(t,y))}{|x-y-2\pi n|^{d+2\lambda}}{\rm d}x {\rm d}y {\rm d}t\\
	&=\sum_{i=1}^{d}\int_0^\tau \int_{\mathbb{T}^d}\int_{(0,2\pi)^d}\vr_\varepsilon(t,x)\vr_\varepsilon(t,y)\frac{(\vu_\varepsilon(t,x)-\vu_\varepsilon(t,y))(\varphi(t,x)-\varphi(t,y))}{|x-y-2\pi e_i|^{d+2\lambda}}{\rm d}x {\rm d}y {\rm d}t\\
	&=\sum_{i=1}^{d}\int_0^\tau \int_{\mathbb{T}^d}\int_{\{(0,2\pi)^d-2\pi e_i\}}\vr_\varepsilon(t,x+2\pi e_i)\vr_\varepsilon(t,y)\frac{(\vu_\varepsilon(t,x+2\pi e_i)-\vu_\varepsilon(t,y))(\varphi(t,x+2\pi e_i)-\varphi(t,y))}{|x-y|^{d+2\lambda}}{\rm d}x {\rm d}y {\rm d}t
	\\
	 &= \int_{0}^\tau \int_{\mathbb{T}^d}\int_{(-2\pi,\,0)^d}\vr_\varepsilon(t,x)\vr_\varepsilon(t,y)\frac{(\vu_\varepsilon(t,x)-\vu_\varepsilon(t,y))(\varphi(t,x)-\varphi(t,y))}{|x-y|^{d+2\lambda}}{\rm d}x {\rm d}y {\rm d}t.
\end{align*}
Now we can conclude that 
\begin{align*}
	I_\varepsilon^{0,1}=\frac{1}{2}\int_0^t\int_{\mathbb{T}^d}\int_{(-2\pi, \,2\pi)^d}\vr_\varepsilon(t,x)\vr_\varepsilon(t,y)\frac{(\vu_\varepsilon(t,x)-\vu_\varepsilon(t,y))(\varphi(t,x)-\varphi(t,y))}{|x-y|^{d+2\lambda}}{\rm d}x {\rm d}y {\rm d}t.
\end{align*} 
We divide the term $I_\varepsilon^{0,1}$ in two parts as 
\begin{align}\label{today0}
	I_\varepsilon^{0,1}=I_{\varepsilon}^{0,m}+ I_\varepsilon^{m,1},
\end{align} 
where
\begin{align*}
	I_\varepsilon^{0,m}&=\frac{1}{2}\int_0^t\int_{\{(x,y)\in(-2\pi,2\pi)^d\times\mathbb{T}^d:|x-y|\textless m\}}\vr_\varepsilon(t,x)\vr_\varepsilon(t,y)\frac{(\vu_\varepsilon(t,x)-\vu_\varepsilon(t,y))(\varphi(t,x)-\varphi(t,y))}{|x-y|^{d+2\lambda}}{\rm d}x {\rm d}y {\rm d}t,\\
	I_\varepsilon^{m,1}&=\frac{1}{2}\int_0^t\int_{\{(x,y)\in(-2\pi,2\pi)^d\times\mathbb{T}^d:|x-y|\ge m\}}\vr_\varepsilon(t,x)\vr_\varepsilon(t,y)\frac{(\vu_\varepsilon(t,x)-\vu_\varepsilon(t,y))(\varphi(t,x)-\varphi(t,y))}{|x-y|^{d+2\lambda}}{\rm d}x {\rm d}y {\rm d}t.
\end{align*}
 By make use of the energy bounds \eqref{energy bounds}, we obtain
\begin{align*}
	|I_{\varepsilon}^{0,m}|=\,&\bigg(\int_0^T\int_{\{(x,y)\in\{(-2\pi,2\pi)^d\times\mathbb{T}^d:|x-y|\textless m\}}\vr_{\varepsilon}(t,x)\vr_{\varepsilon}(t,y)\frac{|\vu_{\varepsilon}(t,x)-\vu_{\varepsilon}(t,y)|^2}{|x-y|^{d+2\lambda}}\dx\dy\dt\bigg)^{1/2}\\&\qquad\bigg(\int_0^T\int_{\{(x,y)\in(-2\pi,2\pi)^d\times\mathbb{T}^d:|x-y|\textless m\}}\vr_{\varepsilon}(t,x)\vr_{\varepsilon}(t,y)\frac{|{\mathbf{ \varphi}}(t,x)-{\mathbf{ \varphi}}(t,y)|^2}{|x-y|^{d+2\lambda}}\dx\dy\dt\bigg)^{1/2}\\&\le\,C\,\bigg(\int_0^T\int_{\{(x,y)\in(-2\pi,2\pi)^d\times\mathbb{T}^d:|x-y|\textless m\}}\vr_{\varepsilon}(t,x)\vr_\varepsilon(t,y)\frac{|{\mathbf{ \varphi}}(t,x)-{\mathbf{ \varphi}}(t,y)|^2}{|x-y|^{d+2\lambda}}\dx\dy\dt\bigg)^{1/2}.
	\end{align*}
We use H\"older inequality to conclude that
\begin{align*}
|I_{\varepsilon}^{0,m}|^2&\le C\,\bigg(\int_0^T\int_{\{(x,y)\in(-2\pi,2\pi)^d\times\mathbb{T}^d:|x-y|\textless m\}}\vr_\varepsilon^2(t,y)\frac{|{\mathbf{ \varphi}}(t,x)-{\mathbf{ \varphi}}(t,y)|^2}{|x-y|^{d+2\lambda}}\dx\dy\dt\bigg)^{1/2}\\&\qquad\qquad\qquad \bigg(\int_0^T\int_{\{(x,y)\in(-2\pi,2\pi)^d\times\mathbb{T}^d:|x-y|\textless m\}}\vr_{\varepsilon}^2(t,x)\frac{|{\mathbf{ \varphi}}(t,x)-{\mathbf{ \varphi}}(t,y)|^2}{|x-y|^{d+2\lambda}}\dx\dy\dt\bigg)^{1/2}
\\&\le C\, \big(\|\vr_\varepsilon\|_{L^\infty([0,T];L^2(\mathbb{T}^d))}m^{1-\lambda}\|\nabla\varphi\|_{L^\infty([0,T]\times\mathbb{T}^d)}\big)\big( \|\vr_\varepsilon\|_{L^\infty([0,T];L^2(\mathbb{T}^d))} m^{1-\lambda}\|\nabla\varphi\|_{L^\infty([0,T]\times\mathbb{T}^d)}\big) \\&\le C\,m^{2-2\lambda}\|\nabla {\mathbf{ \varphi}}\|_{L^\infty([0,T]\times\mathbb{T}^d)}^2\|\vr_{\varepsilon}\|_{L^\infty([0,T];L^2(\mathbb{T}^d))}^2.
\end{align*}
It implies with the help of the energy bounds \eqref{energy bounds} that
\begin{align*}
	|I_{\varepsilon}^{0,m}|\le C\,m^{1-\lambda}\|\nabla {\mathbf{ \varphi}}\|_{L^\infty([0,T]\times\mathbb{T}^d)}.
\end{align*}
It gives that 
\begin{align}\label{today1}
	\lim_{m\to 0}\limsup_{\varepsilon\to 0}I_{\varepsilon}^{0,m}\to 0.
\end{align}
On other hand, for each $m>0$, with the help of \eqref{weak-convergence}, we obtain 
\begin{align}\label{today2}
	\limsup_{\varepsilon\to 0} I_{\varepsilon}^{m,1}=\frac{1}{2}\int_0^{\tau}\int_{\{(x,y)\in(-2\pi,2\pi)^d\times\mathbb{T}^d:|x-y|\ge m\}}\left<\nu_{t,x,y}^2;ss'(\vc{v}-\vc{v'})\right>\frac{\big(\varphi(t,x)-\varphi(t,y)\big)}{|x-y|^{d+2\lambda}}{\rm d}x{\rm d}y {\rm d}t.
\end{align}
By using Lebesgue dominant convergence theorem and \eqref{today0}-\eqref{today2}, we conclude that
\begin{align*}
\limsup_{\varepsilon\to 0} I_{\varepsilon}^{0,1}&=	\lim_{m\to 0}\limsup_{\varepsilon\to 0}I_{\varepsilon}^{0,m} + \lim_{m\to 0}\limsup_{\varepsilon\to 0}I_{\varepsilon}^{m,1}\notag\\&=\frac{1}{2}\int_0^\tau \int_{\mathbb{T}^d}\int_{(-2\pi,2\pi)^d}\left<\nu_{t,x,y}^2;ss'(\vc{v}-\vc{v}')\right>\frac{\big(\varphi(t,x)-\varphi(t,y)\big)}{|x-y|^{d+2\lambda}}\dx\dy\dt\notag\\&=\frac{1}{2}\int_0^\tau \int_{\mathbb{T}^d}\int_{\mathbb{T}^d}\left<\nu_{t,x,y}^2;ss'(\vc{v}-\vc{v}')\right>\frac{\big(\varphi(t,x)-\varphi(t,y)\big)}{|x-y|^{d+2\lambda}}\dx\dy\dt\notag\\&\qquad+\sum_{|n|=1,n\in\mathbb{Z}^d}\frac{1}{2}\int_0^\tau \int_{\mathbb{T}^d}\int_{\mathbb{T}^d}\left<\nu_{t,x,y}^2;ss'(\vc{v}-\vc{v}')\right>\frac{\big(\varphi(t,x)-\varphi(t,y)\big)}{|x-y-2\pi n|^{d+2\lambda}}\dx\dy\dt\notag\\&=\frac{1}{2} \int_0^\tau\int_{\mathbb{T}^d\times\mathbb{T}^d}{\left< \nu_{t,x,y}^2; ss'(\vc{v}-\vc{v}')\right>\big)\big({\mathbf{ \varphi}}(t,x)-{\mathbf{ \varphi}}(t, y)\big)}\kappa_\lambda^{0,1}(x-y)\dx\dy\dt.
\end{align*}
Therefore, we conclude that 
\begin{align}\label{today3}
	\lim_{\varepsilon\to 0}I_\varepsilon^{0,1}=\frac{1}{2} \int_0^\tau\int_{\mathbb{T}^d\times\mathbb{T}^d}{\left< \nu_{t,x,y}^2; ss'(\vc{v}-\vc{v}')\right>\big)\big({\mathbf{ \varphi}}(t,x)-{\mathbf{ \varphi}}(t, y)\big)}\kappa_\lambda^{0,1}(x-y)\dx\dy\dt.
\end{align}
\textbf{Step 2(b)} In this sub-step, we pass the limit in the term $I_\varepsilon^{2,\mathtt{M}}$ as follows: The weak convergence \eqref{weak-convergence} implies that
\begin{align}\label{1st}
	\lim_{\varepsilon\to 0}\int_{[0,\tau]\times\mathbb{T}^d}& {\big(\vr_\epsilon(t,x)\vr_\epsilon(t,y)(\vc{u}_{\varepsilon}(t,x)-\vc{u}_\varepsilon(t,y))\big)\big({\mathbf{ \varphi}}(t,x)-{\mathbf{ \varphi}}(t,y)\big)}\kappa_\lambda^{2,\mathtt{M}}(x-y)\dx\dy\dt\notag\\&=\int_{[0,\tau]\times\mathbb{T}^d} {\left< \nu_{t,x,y}^2; ss'(\vc{v}-\vc{v}')\right>\big({\mathbf{ \varphi}}(t,x)-{\mathbf{ \varphi}}(t,y)\big)}\kappa_\lambda^{2,\mathtt{M}}(x-y)\dx\dy\dt.
\end{align}
With the help of \eqref{1st} and Lebesgue dominant convergence theorem, we have
\begin{align}\label{today8}
	\lim_{\mathtt{M}\to \infty}\lim_{\varepsilon\to 0}I_{\varepsilon}^{2,\mathtt{M}}=\frac{1}{2} \int_0^\tau\int_{\mathbb{T}^d\times\mathbb{T}^d}{\left< \nu_{t,x,y}^2; ss'(\vc{v}-\vc{v}')\right>\big)\big({\mathbf{ \varphi}}(t,x)-{\mathbf{ \varphi}}(t, y)\big)}\big(\kappa_\lambda(x-y)-\kappa_\lambda^{0,1}(x-y)\big)\dx\dy\dt.
\end{align}
\textbf{Step 2(c)} In this sub-step, we estimate the term $I_\varepsilon^{\mathtt{M},\infty}$ as follows: First note that for all $x,y\in(0,2\pi)^d$,
\begin{align}\label{today4}
	|\kappa_\lambda^{\mathtt{M},\infty}(x-y)|\le \sum_{|n|>\mathtt{M}, n\in\mathbb{Z}^d} (2\pi)^{-d-2\lambda}(|n|-1)^{-d-2\lambda}.
\end{align}
$\sum_{|n|>1,n\in\mathbb{Z}^d}(|n|-1)^{-d-2\lambda}$ is a uniformly convergent series, therefore we obtain
\begin{align}\label{today5}
	\lim_{\mathtt{M}\to\infty}\sum_{|n|>\mathtt{M}, n\in\mathbb{Z}^d} (2\pi)^{-d-2\lambda}(|n|-1)^{-d-2\lambda} =0. 
\end{align} 
Again by making use of energy bound \eqref{Energy estimate}, and \eqref{today4}, we have
\begin{align*}
	|I_{\varepsilon}^{\mathtt{M},\infty}|\le\,&\bigg(\int_0^T\int_{\mathbb{T}^d\times \mathtt{T}^d}\vr_{\varepsilon}(t,x)\vr_{\varepsilon}(t,y){|\vu_{\varepsilon}(t,x)-\vu_{\varepsilon}(t,y)|^2}\kappa_\lambda^{\mathtt{M},\infty}(x-y)\dx\dy\dt\bigg)^{1/2}\\&\qquad\bigg(\int_0^T\int_{\mathbb{T}^d\times\mathbb{T}^d}\vr_{\varepsilon}(t,x)\vr_{\varepsilon}(t,y){|{\mathbf{ \varphi}}(t,x)-{\mathbf{ \varphi}}(t,y)|^2}\kappa_\lambda^{\mathtt{M},\infty}(x-y)\dx\dy\dt\bigg)^{1/2}\\&\le\,C\,\bigg(\int_0^T\int_{\mathbb{T}^d\times\mathbb{T}^d}\vr_{\varepsilon}(t,x)\vr_\varepsilon(t,y){|{\mathbf{ \varphi}}(t,x)-{\mathbf{ \varphi}}(t,y)|^2}\kappa_\lambda^{\mathtt{M},\infty}(x-y)\dx\dy\dt\bigg)^{1/2}\\&\le\,C\,\|{\mathbf{ \varphi}}_{\varepsilon}\|_{L^\infty([0,T]\times\mathbb{T}^d)}\|\vr_\varepsilon\|_{L^\infty([0,T];L^2(\mathbb{T}^d))}\bigg( \sum_{|n|>\mathtt{M}, n\in\mathbb{Z}^d} (2\pi)^{-d-2\lambda}(|n|-1)^{-d-2\lambda}\bigg)^{1/2}.
\end{align*}
It shows with the help of \eqref{today5} that
\begin{align}\label{today6}
	\lim_{\mathtt{M}\to\infty}\lim_{\varepsilon\to 0} I_{\varepsilon}^{\mathtt{M},\infty}=0.
\end{align}
By utilizing \eqref{today7},\eqref{today3},\eqref{today8} and \eqref{today6}, we obtain
\begin{align}\label{today15}
	\lim_{\varepsilon\to 0} I_{\varepsilon}=\frac{1}{2} \int_0^\tau\int_{\mathbb{T}^d\times\mathbb{T}^d} {\left< \nu_{t,x,y}^2; ss'(\vc{v}-\vc{v}')\right>\big({\mathbf{ \varphi}}(t,x)-{\mathbf{ \varphi}}(t,y)\big)}{\kappa_\lambda(x-y)}\dx\dy\dt.
\end{align}
\subsection*{Energy inequality.} From the energy estimate \eqref{Energy estimate}, we have
\begin{align*}
	\|\vr_{\varepsilon}(t,x) \vr_{\varepsilon}(t,y) {|\vc{u_{\varepsilon}}(t,x) - \vc{u}_{\varepsilon}(t, y)|^2}\kappa_\lambda(x-y)\|_{ L^1([0,T];L^1(\mathbb{T}^d\times\mathbb{T}^d))}\le\,C.
\end{align*}
Therefore weak-*compactness implies that there exists a $\mu_\lambda\in L^\infty_{\rm weak}([0,T];\mathcal{M}^+(\mathbb{T}^d\times \mathbb{T}^d))$ such that (upto a subsequence)
\begin{align}\label{today9}
	&\frac 12 \int_0^\tau \int_{\mathbb{T}^d \times \mathbb{T}^d } \vr_{\varepsilon}(t,x) \vr_{\varepsilon}(t,y) {|\vc{u_{\varepsilon}}(t,x) - \vc{u}_{\varepsilon}(t, y)|^2}\kappa_\lambda(x-y)\dx\dy\dt\notag\\&\to \frac{1}{2}\int_0^\tau \int_{\mathbb{T}^d\times \mathbb{T}^d} \left<\nu_{t,x,y}^2; {ss'|\vc{v}-\vc{v}'|^2}\right>{\kappa_\lambda(x-y)} \dx\dy\dt +\frac{1}{2} \int_0^\tau\int_{\mathbb{T}^d\times\mathbb{T}^d} d\mu_\lambda(x,y,t)\dt.
\end{align}
\subsection*{Compatibility condition.}
For convenience, we define for all $t\in[0,T]$,
\begin{align*}
	\bar{\vc{v}}(t)=\int_{\mathbb{T}^d}\left<\nu_{t,x};\vc{v}\right>{\rm d}x,\qquad \bar{\vU}(t)=\int_{\mathbb{T}^d}\vU(t,x){\rm d}x,
\end{align*}
and
\begin{align*}
	\bar{\vu}_{\varepsilon}(t)=\int_{\mathbb{T}^d}\vu_\varepsilon(t,x){\rm d}x.
\end{align*}
Let $\mathbf{U}\in C([0,T]; C^1(\mathbb{T}^d))$. By the dominance of semi-norm $\|(-\Delta)^{\lambda}(\cdot)\|_{L^2(\mathbb{T}^d)}$ over the norm $\|\cdot\|_{L^2(\mathbb{T}^d)}$ (see \cite[Proposition 1.2]{BenyiOh2013}), we obtain 
\begin{align*}
	&I_\varepsilon=\int_{0}^\tau\|\big(\mathbf{U}(t)-\bar{\vU}(t)\big)-(\vu_\varepsilon(t)-\bar{\vu}_\varepsilon(t)\big)\|_{L^2(\mathbb{T}^d)}^2\,\dt\\&\le\,C_\lambda\,\int_0^\tau\int_{\mathbb{T}^d\times\mathbb{T}^d}{|\mathbf{U}(t, x)-\vu_\varepsilon(t, x)-\mathbf{U}(t,y)+ \vu_\epsilon(t,y)|^2}\kappa_\lambda(x-y) dx dy dt\\&\le\,\frac{C_\lambda}{c_{\vr}^2}\int_0^\tau\int_{\mathbb{T}^d\times\mathbb{T}^d}{\vr_{\varepsilon}(t,x)\vr_\varepsilon(t,y)|\mathbf{U}(t, x)-\vu_\varepsilon(t, x)-\mathbf{U}(t,y)+ \vu_\epsilon(t,y)|^2}\kappa_\lambda(x-y) dx dy dt
	\end{align*}
\begin{align}\label{today10}&
	=\frac{C_\lambda}{c_{\vr}^2}\int_{0}^\tau\int_{\mathbb{T}^d\times\mathbb{T}^d}\vr_{\varepsilon}(t,x)\vr_{\varepsilon}(t,y)|\mathbf{U}(t,x)-\mathbf{U}(t,y)|^2\kappa_\lambda(x-y) \dx\dy \dt\notag  \\&+\frac{C_\lambda}{c_{\vr}^2} \int_0^\tau \int_{\mathbb{T}^d\times\mathbb{T}^d}\vr_{\varepsilon}(t,x)\vr_{\varepsilon}(t,y)|\vu_{\varepsilon}(t,x)-\vu_\varepsilon(t,y)|^2\kappa_\lambda(x-y)\dx\dy\dt\notag\\&-2\frac{C_\lambda}{c_{\vr}^2}\int_0^\tau\int_{{\mathbb{T}^d}\times\mathbb{T}^d}\vr_{\varepsilon}(t,x)\vr_\varepsilon(t,y)(\vu_{\varepsilon}(t,x)-\vu_{\varepsilon}(t,y))(\mathbf{U}(t,x)-\mathbf{U}(t,y))\kappa_\lambda(x-y)\dx\dy\dt\notag\\
	&=:\frac{C_\lambda}{c_{\vr}^2}(I_1^\varepsilon +  I_2^\varepsilon - I_3^\varepsilon).
\end{align}
By employing the weak convergence of the sequence $\vu_\varepsilon$
and the weak lower semicontinuity of the norm, we obtain 
\begin{align}\label{today11}
	\lim_{\varepsilon\to 0}I_\varepsilon&\ge\int_0^\tau \int_{\mathbb{T}^d}\left<\nu_{t,x}; |\big(\vc{v}-\bar{\vc{v}}(t)\big)-\big(\mathbf{U}(t,x)-\bar{\vU}(t)\big)|^2 \right>\dx\dt.
	\end{align}
We use \eqref{weak-convergence0} to obtain
\begin{align}\label{today12}
	\lim_{\varepsilon\to 0}I_\varepsilon^1&= \int_0^\tau\int_{\mathbb{T}^d\times\mathbb{T}^d} \left< \nu_{t,x,y}^2;ss'|\mathbf{U}(t,x)-\mathbf{U}(t,y)|^2\kappa_\lambda(x-y)\right >\dx\dy\dt.
\end{align}
We use \eqref{today9} to get
\begin{align}\label{today13}
	\lim_{\varepsilon\to 0}I_\varepsilon^2&=\int_0^\tau \int_{\mathbb{T}^d\times\mathbb{T}^d}\left<\nu_{t,x,y}^2;ss'|\vc{v}-\vc{v}'|^2\kappa_\lambda(x-y)\right>\dx\dy\dt +\int_0^\tau \int_{\mathbb{T}^d\times\mathbb{T}^d}{\rm d}\mu_\lambda(x,y,t).
	\end{align}
We use \eqref{today15} to conclude that
\begin{align}\label{today14}
	\lim_{\varepsilon\to 0} I_\varepsilon^3&=2\int_{0}^\tau\int_{\mathbb{T}^d\times\mathbb{T}^d} \left <\nu_{t,x,y}^2; ss'(\vc{v}-\vc{v}')(\mathbf{U}(t,x)-\mathbf{U}(t,y))\kappa_\lambda(x-y)\right>\dx \dy\dt.
\end{align}
We add \eqref{today10}-\eqref{today14} to show that for all $\tau\in [0,T]$,
\begin{align}\label{compatibility first}
	&\int_0^\tau\int_{\mathbb{T}^d}\ \left<\nu_{t,x}; |\big(\vc{v}-\bar{\vc{v}}(t)\big)-(\mathbf{U}(t,x)-\bar{\vU}(t))|^2\right>\dx\dt\notag\\&\le\frac{C_\lambda}{c_{\vr}^2}\int_0^\tau\int_{{\mathbb{T}^d}\times\mathbb{T}^d}\left<\nu_{t,x,y}^2; ss'|(\vc{v}-\vc{v}')-(\mathbf{U}(t,x)-\mathbf{U}(t,y))|^2 \kappa_\lambda(x-y)\right>\dx\dy\dt\notag\\&\qquad +\frac{C_\lambda}{c_{\vr}^2}\int_0^\tau\int_{\mathbb{T}^d\times\mathbb{T}^d}{\rm d}\mu_\lambda(x,y,t).
\end{align}
\noindent
\subsection*{Support of the binary tensor Young measure $\nu$.}
We can prove that 
\begin{align}\label{supp11}
	\mbox{supp}(\nu_{t,x,y}^2)\subset \{(s,\vc{v}, s',\vc{v}'): s,s'\ge c_\vr>0\,\&\, \vc{v},\vc{v}'\in\mathbb{R}^d\}:=A_{t,x,y}^2,
	\end{align}
and 
\begin{align}\label{supp22}
	\mbox{supp}(\nu_{t,x})\subset \{(s,\vc{v}): s\ge c_\vr>0\,\&\, \vc{v}\in\mathbb{R}^d\}:=A_{t,x},
\end{align}
for a.e. $(t,x,y)\in[0,T]\times\mathbb{T}^d\times\mathbb{T}^d.$ Indeed, since $\vr_{\varepsilon} \ge c_\vr>0$, uniformly in $\varepsilon$. Then, we obtain that for every $(s, \vc{v}, s', \vc{v}')\notin A_{t,x,y}^2$, there exists an open set $O_{s,\vc{v}, s', \vc{v}'}$ containing $(s, \vc{v}, s', \vc{v}')$ such that
\begin{align*}
	\nu_{t,x,y}^{2, \varepsilon}\big({O_{s,\vc{v}, s', \vc{v}'}}\big)=0\qquad\forall\,\varepsilon>0.
\end{align*}
It implies that 
\begin{align*}
	\nu_{t,x,y}^2\big({O_{s,\vc{v}, s', \vc{v}'}}\big)\le\,\liminf_{\varepsilon\to 0}\nu_{t,x,y}^{2,\varepsilon}\big({O_{s,\vc{v}, s', \vc{v}'}}\big)=0.
	\end{align*}
It proves \eqref{supp11}.
Similarly, we can establish \eqref{supp22} for $\nu_{t,x}$.
\subsection*{Dissipation defect.}
By the help of weak-*compactness, we can prove that there exist $r^M, r^E\in L^\infty_{\rm weak}([0,T];\mathcal{M}^+(\mathbb{T}^d))$ 
such that
\begin{align*}
	&\frac{1}{2} \vr_{\varepsilon} |\vc{u}_{\varepsilon}|^2 + P(\vr_{\varepsilon})\to  \left< \nu_{t,x}; \frac{1}{2} {s|\vc{v}|^2} + P(s) \right>+ r^E \,\text{weak-* in}\, L^\infty_{\rm weak}([0,T];\mathcal{M}^+(\mathbb{T}^d)),\\
	&\vr_{\varepsilon}\vu_\varepsilon\otimes\vu_\varepsilon\to \left< \nu_{t,x}; {s\vc{v} \otimes \vc{v}}\right> +r^M\, \text{weak-* in}\, L^\infty_{\rm weak}([0,T];\mathcal{M}(\mathbb{T}^d)).
\end{align*}
Now we can set dissipation defect as a.e. $\tau\in[0,T]$
\begin{align*}
	\mathcal {D}(\tau)=\int_{\mathbb{T}^d}{\rm d}r^E(x,\tau)+\frac{1}{2}\int_{0}^\tau \int_{\mathbb{T}^d\times\mathbb{T}^d}{\rm d}\mu_\lambda(x,y,t).
\end{align*}
It is clear from \cite[Lemma 2.1]{Emil_NSE} that defect $\mathcal{D}$ control the contraction measure $r^M$, i.e., \eqref{concentrate defect 2} holds. In our case, $r^C$ is zero. Boundnedess property \eqref{boundedness} is a direct consequence of the energy bounds \eqref{energy bounds}. After passing to the limit in the weak formulations \eqref{viscouse continuity equation}-\eqref{viscous momentum equation} and the energy inequality \eqref{viscouse energy inequality}, we can prove existence of a {\em {\em dissipative measure-valued} } solution to \eqref{NS} in the sense of Definition \ref{DD1}.
\begin{Remark} It is clear from the proof of Theorem \ref{first result} that the energy estimate \eqref{Energy estimate} alone is sufficient to obtain a \emph{dissipative measure-valued} solution in the sense of Definition \ref{DD1}, except the compatibility condition \eqref{compatibility second} and the support of a binary tensor Young measure \eqref{supp1}-\eqref{supp2}. To prove these latter properties, we utilized the uniform lower bound on the densities $\vr_{\varepsilon}$.
\end{Remark}
\section{Proof of Theorem \ref{second result}}\label{Section 4}
\subsection{Relative energy inequality}
\noindent
Note that relative energy functional plays a pivotal role for the comparison of a measure valued solution and a smooth solution.
In our context, we work with the following relative energy functional:
\[
\begin{split}
	&\mathcal{E}_{mv} \left( s, \vu \ \Big| r, \vc{U} \right)(\tau)  = \intO{ \left< \nu_{\tau,x};  \frac{1}{2} s |\vc{v} - \vc{U}(\tau,x)|^2 + P(s) - P'(r(\tau,x)) (s - r(\tau,x)) - P(r(\tau,x)) \right> } \\
	&= \intO{ \left< \nu_{\tau,x}; \frac{1}{2} s |\vc{v}|^2 + P(s) \right> }  - \intO{ \left< \nu_{\tau,x}; s \vc{v} \right> \cdot \vc{U}(\tau,x) } + \intO{ \frac{1}{2}
		\left< \nu_{\tau,x} ; s \right> |\vc{U}(\tau,x)|^2 }\\
	&- \intO{ \left< \nu_{\tau,x} ; s \right> P'(r(\tau,x)) } + \intO{ p(r(\tau,x)) }.
\end{split}
\]
Notice also that for any $\delta>0$ small enough, 
\begin{align}
	\label{imp}
	P(s) - P'(r) (s - r) - P(r) \ge c(\delta) 
	\begin{cases}
		(s -r)^2, \, \text{if}\,\, \delta <r< \delta^{-1},\,s\in[\delta/4,\delta^{-1}] & \\
		(1 + s^\gamma), \,\,\text{if}\,\,\delta<r<\delta^{-1},\, s\notin [\delta/2, 2\delta].&
	\end{cases}
\end{align}
As usual, our main aim is to make use of \ref{dmvEI} and the field equations (\ref{dmvB1}), (\ref{dmvB2}) to deal with the various integral expressions of the relative energy functional. In what follows, we do it in several steps. Firstly, by making use of the equation of continuity (\ref{dmvB1}) with test function $\frac{1}{2}|\vc{U}|^2$, we get
\begin{equation} \label{RE1}
	\begin{split}
		&\intO{ \frac{1}{2} \left< \nu_{\tau,x}; s \right> |\vc{U}|^2(\tau, \cdot) } - \intO{ \frac{1}{2} \left< (\nu_{0})_x; s \right> |\vc{U}|^2(0, \cdot) } \\
		& = \int_0^\tau \intO{ \left[ \left< \nu_{t,x}; s \right> \vc{U}(t,x) \cdot \partial_t \vc{U}(t,x) + \left< \nu_{t,x}; s \vc{v} \right> \cdot \vc{U}(t,x)
			\cdot \Grad \vc{U}(t,x) \right] }\dt \\&\qquad+ \int_0^\tau \intO{ \left< r^C; \frac{1}{2} \Grad |\vc{U}(t,x)|^2 \right> } \dt,
	\end{split}
\end{equation}
provided $\vc{U} \in C^1([0,T] \times {\mathbb{T}^d}; \R^d)$. Again, using $P'(r)$ as a test function in the equation of continuity (\ref{dmvB1}), we get  
\begin{equation} \label{RE2}
	\begin{split}
		& \intO{ \left< \nu_{\tau,x} ;s \right>  P'(r) (\tau, \cdot) }  -  \intO{ \left< (\nu_{0})_x ;s \right>  P'(r) (0, \cdot) }\\
		&= \int_0^\tau \intO{\left[ \left< \nu_{t,x}; s \right> P''(r(t,x))  \partial_t r(t,x) + \left< \nu_{t,x}; s \vc{v} \right> \cdot P''(r(t,x))
			\cdot \Grad r(t,x) \right] }\dt\\&\qquad + \int_0^\tau \intO{ \left< r^C ; P'(r(t,x))  \right> }\dt\\
		&= \int_0^\tau \intO{ \left[ \left< \nu_{t,x}; s \right> \frac{p'(r(t,x))}{r(t,x)}  \cdot \partial_t r(t,x) + \left< \nu_{t,x}; s \vc{v} \right> \frac{p'(r(t,x))}{r(t,x)} \cdot
			\Grad r(t,x) \right] } \dt \\&\qquad+ \int_0^\tau \intO{ \left< r^C ; \Grad P'(r(t,x))  \right> }\dt
	\end{split}
\end{equation} \color{black}
provided $r > 0$ and $r \in C^1([0,T] \times {\mathbb{T}^d})$, and $P$ is twice continuously differentiable in $(0, \infty)$.
Secondly, making use of the momentum equation (\ref{dmvB2}), we have
\begin{equation} \label{RE3}
	\begin{split}
		&\intO{ \left< \nu_{\tau,x}; s \vc{v} \right> \cdot \vc{U} (\tau, \cdot) } - \intO{ \left< (\nu_{0})_x ; s \vc{v} \right> \cdot \vc{U} (0, \cdot) } \\
		&= \int_0^\tau \intO{ \left< \nu_{t,x}; s \vc{v} \right> \cdot \partial_t \vc{U}(t,x) }\dt + \int_0^\tau \int_{{\mathbb{T}^d}} \left[ \left< \nu_{t,x}; s \vc{v} \otimes \vc{v} \right> : \Grad \vc{U}(t,x) + \left< \nu_{t,x}; p(s) \right> \Div \vc{U}(t,x) \right]\dx\dt\\
		& - \frac{1}{2}\int_0^\tau \int_{\mathbb{T}^d \times \mathbb{T}^d} \left< \nu_{t,x,y}^2 ; s s' {(\vc{v} - \vc{v}') (\vc{U}(t,x)-\vc{U}(t,y))} \kappa_\lambda(x-y)\right> \dx\dy\dt+ \int_0^\tau \left< r^M; \Grad \vc{U}(t,x) \right>  \ \dt,
	\end{split}
\end{equation}
for any smooth function $\vc{U} \in C^1([0,T] \times {\mathbb{T}^d}; \R^d)$. In view of the three expressions \eqref{RE1}--\eqref{RE3}, we can recast the relative energy functional as follows:
\begin{equation} \label{RE5}
	\begin{split}
		\mathcal{E}_{mv} \left( s, \vu \ \Big| r , \vc{U} \right)    &+\frac{1}{2}
		\int_0^\tau \int_{\mathbb{T}^d \times \mathbb{T}^d} \left< \nu_{t,x}^2; s s' {(\vc{v} - \vc{v}')} \left((\vc{v} - \vc{v}')- (\vc{U}(t,x)-\vc{U}(t,y))\kappa_\lambda(x-y) \right)\right> \dx\dy\dt \\
		+ \mathcal{D}(\tau)&
		\leq \intO{ \left< (\nu_{0})_x (s - r_0(x)) - P(r_0(x)) \right> }- \int_0^\tau \intO{ \left< \nu_{t,x}, s \vc{v} \right> \cdot \partial_t \vc{U}(t,x) }\dt\\
		& - \int_0^\tau \int_{{\mathbb{T}^d}} \left[ \left< \nu_{t,x};  s \vc{v} \otimes \vc{v} \right> : \Grad \vc{U}(t,x) + \left< \nu_{t,x}; p(s) \right> \Div \vc{U}(t,x) \right] \dx\dt\\
		& + \int_0^\tau \intO{ \left[ \left< \nu_{t,x}; s \right> \vc{U}(t,x)  \cdot \partial_t \vc{U}(t,x) + \left< \nu_{t,x}; s \vc{v} \right> \cdot \vc{U}(t,x)
			\cdot \Grad \vc{U}(t,x) \right] }\dt\\
		&+ \int_0^\tau \intO{ \left[ \left< \nu_{t,x} ; \left(1 - \frac{s}{r(t,x)}  \right) \right> p'(r(t,x)) \partial_t r - \left< \nu_{t,x}; s \vc{v} \right> \cdot \frac{p'(r(t,x))}{r(t,x)}
			\Grad r(t,x) \right] }\dt\\
		&+ \int_0^\tau  \left< r^C; \frac{1}{2} \Grad |\vc{U}(t,x)|^2  - \Grad P'(r(t,x)) \right>  \  \dt - \int_0^\tau  \left< r^M ; \Grad \vc{U}(t,x) \right>\dt.
	\end{split}
\end{equation}
To deal with the last two terms of \eqref{RE5}, we recall Definition \ref{DD1} to conclude
\[
\begin{split}\label{TEST}
	&\left| \int_0^\tau  \left< r^C; \frac{1}{2} \Grad |\vc{U}|^2  - \Grad P'(r) \right>  \  \dt - \int_0^\tau  \left< r^M ; \Grad \vc{U} \right>  \dt \right| \\
	&\leq C \left( \left\| \Grad \vc{U} \right\|_{C([0, T] \times {\mathbb{T}^d}; \mathbb{R}^{d \times d})} +
	\left\| r \right\|_{C([0,T] \times {\mathbb{T}^d}) } + \left\| \Grad r \right\|_{C([0, T] \times {\mathbb{T}^d}; \mathbb{R}^{d})}    \right) \int_0^\tau (\chi(t)+\xi(t)) \mathcal{D}(t) \ \dt.
\end{split}
\]
Next, we make use of continuity equation satisfied by $(r, \vc{U})$, i.e.,
\bFormula{RE6}
\partial_t r + \Div (r \vc{U} ) = 0
\eF
to rewrite \eqref{RE5} as
\begin{equation} \label{RE7}
	\begin{split}
		\mathcal{E}_{mv} \left( s, \vu \ \Big| r , \vc{U} \right)    &+\frac{1}{2}
		\int_0^\tau \int_{\mathbb{T}^d \times \mathbb{T}^d} \left< \nu_{t,x,y}^2 ; s s'{(\vc{v} - \vc{v}')} \left((\vc{v} - \vc{v}')- (\vc{U}(t,x)-\vc{U}(t,y)) \right)\kappa_\lambda(x-y)\right>\dx\dy\dt  \\
		&+ \mathcal{D}(\tau) 
		\leq \intO{ \left< (\nu_{0})_x;  \frac{1}{2} s |\vc{v} - \vc{U}_0(x)|^2 + P(s) - P'(r_0(x)) (s - r_0(x)) - P(r_0(x)) \right> } \\
		& - \int_0^\tau \intO{ \left< \nu_{t,x}, s \vc{v} \right> \cdot \partial_t \vc{U}(t,x) }\dt - \int_0^\tau \int_{{\mathbb{T}^d}}  \left< \nu_{t,x};  s \vc{v} \otimes \vc{v} \right> : \Grad \vc{U}(t,x)  \dx\dt\\
		& + \int_0^\tau \intO{ \left[ \left< \nu_{t,x}; s \right> \vc{U}(t,x)  \cdot \partial_t \vc{U}(t,x) + \left< \nu_{t,x}; s \vc{v} \right> \cdot \vc{U}(t,x)
			\cdot \Grad \vc{U}(t,x) \right] }\dt\\
		& + \int_0^\tau \intO{  \left< \nu_{t,x}; s \vc{U}(t,x) - s \vc{v} \right> \cdot \frac{p'(r(t,x))}{r(t,x)} \Grad r(t,x)  } \ \dt \\
		&- \int_0^\tau \intO{  \left< \nu_{t,x}; p(s) - p'(r(t,x))(s -r(t,x)) - p(r(t,x))  \right> \Div \vc{U}(t,x) } \dt\\
		&+ C\, \int_0^\tau  (\chi(t)+\xi(t)) \mathcal{D}(t)\dt.
	\end{split}
\end{equation}
To simplify further, we make use of the momentum equation satisfied by $(r, \vc{U})$, i.e., 
\begin{align}\label{expression}
	\partial_t \vc{U}(t,x) + \vc{U}(t,x) \cdot \Grad \vc{U}(t,x) + \frac{1}{r} \Grad p(r(t,x)) = \int_{\mathbb{T}^d} {\big(\vc{U}(t,y)-\vc{U}(t,x)\big)}\kappa_\lambda(x-y) r(t,y)\,dy,
\end{align}
and rewrite \eqref{RE7} as
\begin{equation} \label{RE8}
	\begin{split}
		&\mathcal{E}_{mv} \left( s, \vu \ \Big| r , \vc{U} \right)+
		\frac{1}{2}\int_0^\tau \int_{\mathbb{T}^d \times \mathbb{T}^d} \left< \nu_{t,x,y}^2 ; s s'{(\vc{v} - \vc{v}')} \left((\vc{v} - \vc{v}')- (\vc{U}(t,x)-\vc{U}(t,y)) \right)\kappa_\lambda(x-y)\right> \dx\dy\dt \\
		&+ \mathcal{D}(\tau) 
		\leq \intO{ \left< (\nu_{0})_x;  \frac{1}{2} s |\vc{v} - \vc{U}_0(x)|^2 + P(s) - P'(r_0) (s - r_0(x)) - P(r_0(x)) \right> } \\
		& + \int_0^\tau \intO{ \left< \nu_{t,x}, s \vc{v}  - s \vc{U}(t,x) \right> \cdot \Grad \vc{U}(t,x) \cdot \vc{U}(t,x) } \ \dt - \int_0^\tau \int_{{\mathbb{T}^d}}  \left< \nu_{t,x};  s \vc{v} \otimes \vc{v} \right> : \Grad \vc{U}(t,x) \dx \dt\\
		& + \int_0^\tau \intO{ \left< \nu_{t,x}; s \vc{v} \right> \cdot
			\Grad \vc{U}(t,x)\cdot  \vc{U}(t,x)}\dt\\
		& + \int_0^\tau \intO{  \left< \nu_{t,x}; s \vc{U}(t,x) - s \vc{v} \right> \cdot \int_{\mathbb{T}^d}{\big(\vc{U}(t,y)-\vc{U}(t,x)\big)}\kappa_\lambda(x-y) r(t,y)\,\dy \, }\dt \\
		&+ \int_0^\tau \intO{  \left< \nu_{t,x}; p(s) - p'(r(t,x))(s -r(t,x)) - p(r(t,x))  \right> \Div \vc{U}(t,x) }\dt\\
		&+ C\, \int_0^\tau  (\chi(t)+\xi(t)) \mathcal{D}(t) \dt.
	\end{split}
\end{equation}
Next, we rewrite 
\begin{align*}
	&s s' {(\vc{v} - \vc{v}')}\left((\vc{v} - \vc{v}')- (\vc{U}(t,x)-\vc{U}(t,y)) \right)\\
	&= \Big(s s'{(\vc{v} - \vc{v}')} -  r(t,x) r(t,y){(\vc{U}(t,x) - \vc{U}(t,y))}\Big)\left((\vc{v} - \vc{v}')- (\vc{U}(t,x)-\vc{U}(t,y)) \right) \\
	&\qquad + r(t,x) r(t,y){(\vc{U}(t,x) - \vc{U}(t,y))} \left((\vc{v} - \vc{v}')- (\vc{U}(t,x)-\vc{U}(t,y)) \right).
\end{align*}
Moreover, we can write
\begin{align}\label{today16}
	&\Big(s s' {(\vc{v} - \vc{v}')} -  r(t,x) r(t,y) {(\vc{U}(t,x) - \vc{U}(t,y))} \Big)\left((\vc{v} - \vc{v}')- (\vc{U}(t,x)-\vc{U}(t,y)) \right)\notag \\
	&= s s' {\left((\vc{v} - \vc{v}')- (\vc{U}(t,x)-\vc{U}(t,y)) \right)} \left((\vc{v} - \vc{v}')- (\vc{U}(t,x)-\vc{U}(t, y)) \right)\notag\\
	&\qquad + (ss'-r(t,x)r(t,y)){(\vc{U}(t,x) - \vc{U}(t,y))}
	\left((\vc{v} - \vc{v}')- (\vc{U}(t,x)-\vc{U}(t,y)) \right).
\end{align}
We also have
\begin{align}\label{today17}
\begin{split}
	&\int_0^\tau \intO{ \left< \nu_{t,x}; \left( s \vc{v} - s \vc{U}(t,x) \right) \right> \cdot \Grad \vc{U}(t,x) \cdot \vc{U}(t,x) } \ \dt
	- \int_0^\tau \int_{{\mathbb{T}^d}} \left< \nu_{t,x}; s \vc{v} \otimes \vc{v} \right> : \Grad \vc{U}(t,x) \dx\dt\\
	&+ \int_0^\tau \intO{ \left< \nu_{t,x}; s \vc{v} \right> \cdot \Grad \vc{U}(t,x) \cdot \vc{U}(t,x)   } \ \dt\\
	&=  \int_0^\tau \intO{ \left< \nu_{t,x}; \left( s \vc{v} - s \vc{U}(t,x) \right) \right> \cdot \Grad \vc{U}(t,x) \cdot \vc{U}(t,x) }\dt+
	\int_0^\tau \int_{{\mathbb{T}^d}} \left< \nu_{t,x}; s \vc{v} \cdot   \Grad \vc{U}(t,x)\cdot(\vc{U}(t,x) - \vc{v} ) \right> \dx\dt\\
	&= \int_0^\tau \int_{{\mathbb{T}^d}} \left< \nu_{t,x}; s (\vc{v} - \vc{U}(t,x)) \cdot \Grad \vc{U}(t,x)\cdot (\vc{U}(t,x) - \vc{v} )  \right>  \dx \dt.
\end{split}
\end{align}
Therefore, (\ref{RE8}) can be written with the help of \eqref{today16}-\eqref{today17} as
\begin{align} 
		&\mathcal{E}_{mv} \left( s, \vu \ \Big| r , \vc{U} \right)(\tau)  +\frac{1}{2}
		\int_0^\tau \int_{\mathbb{T}^d \times \mathbb{T}^d} \left< \nu_{t,x,y}^2; s s' {|(\vc{v} - \vc{v}')- (\vc{U}(t,x)-\vc{U}(t,y))|^2}k_{\lambda}(x-y)\right>\dx\dy\dt \nonumber \\
		&+ \mathcal{D}(\tau) 
		\leq \intO{ \left< (\nu_{0})_x;  \frac{1}{2} s |\vc{v} - \vc{U}_0|^2 + P(s) - P'(r_0(x)) (s - r_0(x)) - P(r_0(x)) \right> } \nonumber \\
		& + \int_0^\tau \int_{{\mathbb{T}^d}} \left< \nu_{t,x}; s (\vc{v} - \vc{U}(t,x))\cdot \Grad \vc{U}(t,x) \cdot (\vc{U}(t,x) - \vc{v} )  \right>  \dx\dt \label{RE9} \\
		& -\frac{1}{2} 
		\int_0^\tau \int_{\mathbb{T}^d \times \mathbb{T}^d} \left< \nu_{t,x,y}^2 ; r(t,x) r(t,y){(\vc{U}(t,x) - \vc{U}(t,y))} \kappa_\lambda(x-y)\left((\vc{v} - \vc{v}')- (\vc{U}(t,x)-\vc{U}(t,y)) \right)\right> \dx\dy\dt \nonumber \\
		& + \int_0^\tau \intO{  \left< \nu_{t,x}; s \vc{U}(t,x) - s \vc{v} \right> \cdot \int_{\mathbb{T}^d} {\vc{U}(t,y)-\vc{U}(t,x)}\kappa_\lambda(x-y) r(y,t)\,\dy\,  } \ \dt \nonumber \\
		&+ \int_0^\tau \intO{  \left< \nu_{t,x} ; p(s) - p'(r(t,x))(s -r(t,x)) - p(r(t,x))  \right> \Div \vc{U}(t,x) } \ \dt \nonumber \\
		& - \int_0^\tau \int_{\mathbb{T}^d \times \mathbb{T}^d} \left< \nu_{t,x,y}^2 ; (ss'-r(t,x)r(t,y)){(\vc{U}(t,x) - \vc{U}(t,y))}k_{\lambda}(x-y)
		\left((\vc{v} - \vc{v}')- (\vc{U}(t,x)-\vc{U}(t, y)) \right)\right> \dx\dy\dt \nonumber \\
		&+ C\, \int_0^\tau  (\chi(t)+\xi(t)) \mathcal{D}(t)  \dt. \nonumber
\end{align}
Our next goal is to deal with the following terms:
\begin{align*}
	&\int_0^\tau \intO{  \left< \nu_{t,x}; s \vc{U}(t,x) - s \vc{v} \right> \cdot \int_{\mathbb{T}^d} {\vc{U}(t,y)-\vc{U}(t,x)}\kappa_\lambda(x-y) r(t,y)\,\dy\,  } \ \dt \\
	& \quad -\frac{1}{2}\int_0^\tau \int_{\mathbb{T}^d \times \mathbb{T}^d} \left< \nu_{t,x,y}^2; r(t,x) r(t,y) {(\vc{U}(t,x) - \vc{U}(t,y))}\kappa_\lambda(x-y) \left((\vc{v} - \vc{v}')- (\vc{U}(t,x)-\vc{U}(t,y)) \right)\right>\dx\dy\dt .
\end{align*}
Notice that, by interchanging the role of coordinates, we have 
\begin{align}\label{today 001}
	&-\frac{1}{2} \int_0^\tau \int_{\mathbb{T}^d \times \mathbb{T}^d} \left< \nu_{t,x,y}^2 ; r(t,x) r(t,y) {(\vc{U}(t,x) - \vc{U}(t,y))}k_{\lambda}(x-y) \left((\vc{v} - \vc{v}')- (\vc{U}(t,x)-\vc{U}(t,y)) \right)\right>\dx\dy\dt\notag\\
	&
	=-\frac{1}{2} \int_0^\tau \int_{\mathbb{T}^d \times \mathbb{T}^d} \left< \nu_{t,x,y}^2 ; r(t,x) r(t,y) {(\vc{U}(t,x) - \vc{U}(t,y))}k_{\lambda}(x-y) (\vc{v}-\vc{U}(t,x)) \right> \,\dy\dx\dt\notag\\
	&\qquad+\frac{1}{2} \int_0^\tau \int_{\mathbb{T}^d \times \mathbb{T}^d} \left< \nu_{t,x,y}^2 ; r(t,x) r(t,y) {(\vc{U}(t,x) - \vc{U}(t,y))}k_{\lambda}(x-y) (\vc{v}'-\vc{U}(t,y)) \right> \dy\,\dx\dt\notag \\&
	=-\frac{1}{2} \int_0^\tau \int_{\mathbb{T}^d \times \mathbb{T}^d} \left< \nu_{t,x,y}^2 ; r(t,x) r(t,y) {(\vc{U}(t,x) - \vc{U}(t,y))}k_{\lambda}(x-y) (\vc{v}-\vc{U}(t,x)) \right>\dy\dx\dt\notag\\
	&\qquad-\frac{1}{2} \int_0^\tau \int_{\mathbb{T}^d \times \mathbb{T}^d} \left< \nu_{t,y,x}^2 ; r(t,x) r(t,y) {(\vc{U}(t,x) - \vc{U}(t,y))}k_{\lambda}(x-y) (\vc{v}'-\vc{U}(t,x)) \right> \dy\dx\dt.\notag
	\end{align}
By using the symmetry property of binary tensor Young measure $\nu$ (see Remark \ref{remark1}) and \eqref{L^2 boundedness}, we conclude that
\begin{align}&-\frac{1}{2} \int_0^\tau \int_{\mathbb{T}^d \times \mathbb{T}^d} \left< \nu_{t,x,y}^2 ; r(t,x) r(t,y) {(\vc{U}(t,x) - \vc{U}(t,y))}k_{\lambda}(x-y) \left((\vc{v} - \vc{v}')- (\vc{U}(t,x)-\vc{U}(t,y)) \right)\right>\dx\dy\dt\notag\\&=-\frac{1}{2} \int_0^\tau \int_{\mathbb{T}^d \times \mathbb{T}^d} \left< \nu_{t,x,y}^2 ; r(t,x) r(t,y) {(\vc{U}(t,x) - \vc{U}(t,y))}k_{\lambda}(x-y) (\vc{v}-\vc{U}(t,x)) \right>\dy\dx\dt\notag\\
	&\qquad-\frac{1}{2} \int_0^\tau \int_{\mathbb{T}^d \times \mathbb{T}^d} \left< \nu_{t,x,y}^2 ; r(t,x) r(t,y) {(\vc{U}(t,x) - \vc{U}(t,y))}k_{\lambda}(x-y) (\vc{v}-\vc{U}(t,x)) \right> \dy\dx\dt\notag\\
	&=\int_0^\tau \int_{\mathbb{T}^d \times \mathbb{T}^d} \left< \nu_{t,x,y}^2 ; r(t,x) r(t,y) {(\vc{U}(t,y) - \vc{U}(t,x))}k_{\lambda}(x-y) (\vc{v}-\vc{U}(t,x)) \right>\dy\dx\dt.
	\end{align}
Note that, we also have
\begin{align}\label{today 002}
	&\int_{\mathbb{T}^d \times \mathbb{T}^d} \left< \nu_{t,x};  \left( r(t,x) \vc{v} - r(t,x) \vc{U}(t,x) \right) \right> \cdot  {\big(\vc{U}(t,x)-\vc{U}(t,y)\big)}\kappa_\lambda(x-y) r(t,y)\dx \ \dy\notag \\
	&= \int_{\mathbb{T}^d \times \mathbb{T}^d} \left< \nu_{t,x};  \left( r(t,x) \vc{v} - r(t,x) \vc{U}(t,x) \right) \right> \cdot \frac{1}{r(t,x)}{\big(\vc{U}(t,x)-\vc{U}(t,y)\big)}\kappa_\lambda(x-y) r(t,x) r(t,y)\dx \ \dy\notag \\
	& = \int_{\mathbb{T}^d \times \mathbb{T}^d} \left<\nu_{t,x};  r(t,x)\vc{v} - r(t,x)\vc{U}(t,x)\right> \cdot \frac{1}{r(t,x)} {\big(\vc{U}(t,x)-\vc{U}(t,y)\big)}\kappa_\lambda(x-y) r(t,x) r(t,y)\dx \ \dy\notag \\
	& =  \int_{\mathbb{T}^d \times \mathbb{T}^d} \left<\nu_{t,x};( \vc{v} - \vc{U}(t,x)) \right> \cdot {\big(\vc{U}(t,x)-\vc{U}(t,y)\big)}\kappa_\lambda(x-y) r(t,x) r(t,y)\dx \ \dy\notag
	\\
	& =  \frac12 \int_{\mathbb{T}^d \times \mathbb{T}^d}\left<\nu_{t,x,y}^2; r(t,x) r(t,y){(\vc{U}(t,x) - \vc{U}(t,y))}\kappa_\lambda(x-y) \left((\vc{v} - \vc{v}')- (\vc{U}(t,x)-\vc{U}(t,y)) \right)\right> \dx\dy.
\end{align}
Therefore, by using the consistency property of the tensor Young measure $\nu$ (see Remark \ref{remark1}), the identities \eqref{today 001}-\eqref{today 002} give 
\begin{align}\label{today 003}
\begin{split}
	&\int_0^\tau \intO{  \left< \nu_{t,x,y}^2; s \vc{U}(t,x) - s \vc{v} \right> \cdot \int_{\mathbb{T}^d} {\big(\vc{U}(t,y)-\vc{U}(t,x)\big)}\kappa_\lambda(x-y) r(t,y)\dy}\dt \\
	& \quad -\frac{1}{2} \int_0^\tau \int_{\mathbb{T}^d \times \mathbb{T}^d} \left< \nu_{t,x,y}^2 ; r(t,x) r(t,y) {(\vc{U}(t,x) - \vc{U}(t,y))}k_{\lambda}(x-y) \left((\vc{v} - \vc{v}')- \big(\vc{U}(t,x)-\vc{U}(t,y)\big) \right)\right> \dx\dy\dt\\
	& = \int_0^\tau \intO{ \left< \nu_{t,x};  \left( s \vc{U}(t,x) - s \vc{v} + r(t,x) \vc{v} - r(t,x) \vc{U}(t,x) \right) \right> \cdot \int_{\mathbb{T}^d} {\big(\vc{U}(t,y)-\vc{U}(t,x)\big)}\kappa_\lambda(x-y) r(t,y)\,\dy }\dt\\
	& = \int_0^\tau \intO{ \left< \nu_{t,x};  (s - r(t,x)) (\vc{U}(t,x) - \vc{v} ) \right> \cdot \int_{\mathbb{T}^d} {\big(\vc{U}(t,y)-\vc{U}(t,x)\big)} \kappa_\lambda(x-y) r(t,y)\,\dy }\dt.
\end{split}
\end{align}
We use \eqref{RE9} and \eqref{today 003} to conclude that
\begin{equation} \label{RE100}
	\begin{split}
		&\mathcal{E}_{mv} \left( s, \vu \ \Big| r , \vc{U} \right)(\tau)  +\frac{1}{2}
		\int_0^\tau \int_{\mathbb{T}^d \times \mathbb{T}^d} \left< \nu_{t,x,y}^2; s s' {|(\vc{v} - \vc{v}')- (\vc{U}(t,x)-\vc{U}(t,y))|^2}k_{\lambda}(x-y)\right>\dx\dy\dt \\
		&+ \mathcal{D}(\tau) 
		\leq \intO{ \left< (\nu_{0})_x;  \frac{1}{2} s |\vc{v} - \vc{U}_0|^2 + P(s) - P'(r_0(x)) (s - r_0(x)) - P(r_0(x)) \right> } + C\, \int_0^\tau  (\chi(t)+\xi(t)) \mathcal{D}(t)  \dt\\&\qquad\qquad +I_1 +I_2 +I_3 +I_4,
	\end{split}
\end{equation}
where
\begin{align*}
	I_1 &= \int_0^\tau \int_{{\mathbb{T}^d}} \left< \nu_{t,x}; s (\vc{v} - \vc{U}(t,x))\cdot \Grad \vc{U}(t,x) \cdot (\vc{U}(t,x) - \vc{v} )  \right>  \dx\dt,\\
	I_2&=  \int_0^\tau \intO{  \left< \nu_{t,x} ; p(s) - p'(r(t,x))(s -r(t,x)) - p(r(t,x))  \right> \Div \vc{U}(t,x) } \ \dt,\\
	I_3&= \int_0^\tau \intO{ \left< \nu_{t,x};  (s - r(t,x)) (\vc{U}(t,x) - \vc{v} ) \right> \cdot \int_{\mathbb{T}^d} \big({\vc{U}(t,y)-\vc{U}(t,x)}\big)\kappa_\lambda(x-y) r(t,y)\,\dy } \ \dt,\\
	I_4&=-\int_0^\tau\int_{\mathbb{T}^d\times\mathbb{T}^d}\big\langle \nu_{t,x,y}^2; (ss'-r(t,x)r(t,y)){(\vc{U}(t,x) - \vc{U}(t,y))}\kappa_\lambda(x-y)\\&\qquad\qquad\qquad\qquad\qquad\qquad\qquad\qquad
	\left((\vc{v} - \vc{v}')-(\vc{U}(t,x)-\vc{U}(t,y)) \right)\big\rangle\dx\dy\dt.
\end{align*}
We estimate these above terms individually through the following steps:

\noindent
\textbf{Step 1.}
For terms $I_1 $ and $I_2$, we use \eqref{assumption} to conclude that
\begin{align}\label{today101}
	I_1 + I_2&\le \|\nabla_x\vU\|_{L^\infty([0,T]\times\mathbb{T}^d)} \int_0^\tau \mathcal{E}_{mv} \left( s, \vu \ \Big| r , \vc{U} \right)(t){\rm d}t\notag\\&\le\, C\,\int_0^\tau \mathcal{E}_{mv} \left( s, \vu \ \Big| r , \vc{U} \right)(t){\rm d}t.
\end{align}
\textbf{Step 2.}
For the term $I_3$, by using the first compatibility condition \eqref{compatibility first} and a similar argument as used in \cite[page 15, Section 4.3]{Emil_NSE}, we can conclude that for any $\delta_1>0$
	\begin{align}\label{today000}
		I_3&\le\,C_\lambda(\delta_1)\,\int_0^\tau \mathcal{E}_{mv}(s,\vu\big|r,\mathbf{U})(t) \dt+\delta_1 C_\lambda \mathcal{D}(\tau)\notag\\&\qquad\qquad+\delta_1 \int_0^\tau\int_{{\mathbb{T}^d}\times\mathbb{T}^d}\left<\nu_{t,x,y}^2; ss'|(\vc{v}-\vc{v}')-(\mathbf{U}(t,x)-\mathbf{U}(t,y))|^2 \kappa_\lambda(x-y)\right>\dx\dy\dt + \delta_1 \mathcal{D}(\tau).
	\end{align}
Indeed, we write
\begin{align}\label{today0000}
	&\int_0^\tau \intO{ \left< \nu_{t,x};  (s - r(t,x)) (\vc{U}(t,x) - \vc{v} ) \right> \cdot \int_{\mathbb{T}^d} {\big(\vc{U}(t,y)-\vc{U}(t,x)\big)}\kappa_\lambda(x-y) r(t,y)\,\dy } \dt\notag\\
	&=\int_0^\tau \intO{ \left< \nu_{t,x}; \psi(s) (s - r(t,x)) (\vc{U}(t,x) - \vc{v} ) \right> \cdot \int_{\mathbb{T}^d} {\big(\vc{U}(t,y)-\vc{U}(t,x)\big)}\kappa_\lambda(x-y) r(t,y)\,\dy } \dt\notag\\&\qquad+\int_0^\tau \intO{ \left< \nu_{t,x}; \big(1-\psi(s)\big) (s - r(t,x)) (\vc{U}(t,x) - \vc{v} ) \right> \cdot \int_{\mathbb{T}^d} {\big(\vc{U}(t,y)-\vc{U}(t,x)\big)}\kappa_\lambda(x-y) r(t,y)\,\dy } \dt\notag\\&=:I_\psi + I_{(1-\psi)},
\end{align}
where $$\psi\in C_c^\infty(0,\infty),\qquad 0\le\psi\le 1,\qquad\psi(s)=1 \qquad\text{for}\qquad s\in (\inf r,\sup r).$$
We choose $\delta>0$ small enough such that 
\begin{align}\label{support0}
	\mbox{supp}[\psi]\subseteq [\delta/4,\delta^{-1}],\qquad [\delta/2, \delta^{-1}/2]\subseteq [\inf r, \sup r]\subseteq[\delta,\delta^{-1}],\qquad\text{and}\qquad 4\delta< \delta^{-1}.
	\end{align}
Consequently, for the term $I_\psi$, we get 
\begin{align}\label{today0002}
	&\int_0^\tau \intO{ \left< \nu_{t,x}; \psi(s) (s - r(t,x)) (\vc{U}(t,x) - \vc{v} ) \right> \cdot \int_{\mathbb{T}^d} {\big(\vc{U}(t,y)-\vc{U}(t,x)\big)}\kappa_\lambda(x-y) r(t,y)\,\dy } \dt\notag \\&\le\sup_{(x,t)\in \mathbb{T}^d\times[0,T]}\bigg|\int_{\mathbb{T}^d} {\big(\vc{U}(t,y)-\vc{U}(t,x)\big)}\kappa_\lambda(x-y) r(t,y)\,\dy\bigg|\notag\\&\qquad\qquad\qquad\qquad\qquad\int_0^\tau \intO{ \left< \nu_{t,x}; \psi(s) |(s - r(t,x))| |(\vc{U}(t,x) - \vc{v} )| \right> } \dt.
\end{align}
By using \eqref{assumption} and \eqref{expression}, we conclude that there exists a constant $C>0$ such that
\begin{align}\label{today0001}
	\sup_{(x,t)\in \mathbb{T}^d\times[0,T]}\bigg|\int_{\mathbb{T}^d} {\big(\vc{U}(t,y)-\vc{U}(t,x)\big)}\kappa_\lambda(x-y) r(t,y)\,\dy\bigg|\le\,C.
\end{align}
We also have 
\begin{align*}
\int_0^\tau \intO{ \left< \nu_{t,x}; \psi(s) |(s - r(t,x))| |(\vc{U}(t,x) - \vc{v} )| \right> } \dt\le&\,\frac{1}{2}\int_0^\tau \intO{ \left< \nu_{t,x};\frac{\psi(s)^2}{{s}} |(s - r(t,x))|^2  \right> } \dt\notag\\&\,\,+\frac{1}{2}\int_0^\tau \intO{ \left< \nu_{t,x};s|(\vc{U}(t,x) - \vc{v} )|^2  \right> } \dt.
\end{align*}
With the help of \eqref{imp} and \eqref{support0}, it implies that 
\begin{align}\label{today0003}
	\int_0^\tau \intO{ \left< \nu_{t,x}; \psi(s) |(s - r(t,x))| |(\vc{U}(t,x) - \vc{v} )| \right> } \dt&\le\, C\,\int_0^\tau \mathcal{E}_{mv} \left( s, \vu \ \Big| r , \vc{U} \right)(t){\rm d}t.
\end{align}
We add \eqref{today0000}-\eqref{today0003} to obtain
\begin{align}\label{today0004}
	|I_{\psi}|\le C\, \int_0^\tau \mathcal{E}_{mv} \left( s, \vu \ \Big| r , \vc{U} \right)(t){\rm d}t.
\end{align}
For term $I_{(1-\psi)}$, next we write 
\begin{align}\label{today0005}
	I_{(1-\psi)}=I_{\omega_1}+I_{\omega_2},
\end{align}
where
\begin{align*}
	I_{\omega_1}&=\int_0^\tau \intO{ \left< \nu_{t,x}; \omega_1(s)(s - r(t,x)) (\vc{U}(t,x) - \vc{v} ) \right> \cdot \int_{\mathbb{T}^d} {\big(\vc{U}(t,y)-\vc{U}(t,x)\big)}\kappa_\lambda(x-y) r(t,y)\,\dy } \dt,\\
	I_{\omega_2}&=\int_0^\tau \intO{ \left< \nu_{t,x}; \omega_2(s)(s - r(t,x)) (\vc{U}(t,x) - \vc{v} ) \right> \cdot \int_{\mathbb{T}^d} {\big(\vc{U}(t,y)-\vc{U}(t,x)\big)}\kappa_\lambda(x-y) r(t,y)\,\dy } \dt.
\end{align*}
Here 
$$\mbox{supp}[\omega_1]\subset[0,\inf r),\qquad\mbox{supp}[\omega_2]\subset[\sup r,\infty),\qquad\omega_1+ \omega_2=1-\psi.$$
It is clear from \eqref{support0} that 
\begin{align}\label{support}
	\mbox{supp}[\omega_1]\subseteq [0,\delta/2)\qquad\text{and}\qquad \mbox{supp}[\omega_2]\subseteq[\delta^{-1}/2,\infty).
\end{align}
Accordingly, for small $\delta_1>0$, we have
\begin{align*}
	\left<\nu_{t,x},\omega_1(s)|(s-r)||(\vU-\vc{v})|\right>\le &C_\lambda(\delta_1)\left<\nu_{t,x};\omega_1^2(s)(s-r)^2\right> + \frac{\delta_1}{C_\lambda}\left<\nu_{t,x};|\big(\vU(t)-\bar{\vU}(t)\big)-\big(\vc{v}-\bar{\vc{v}}(t)\big)|^2\right>\\&\qquad+\left<\nu_{t,x},\omega_1(s)|(s-r)||(\bar{\vU}(t)-\bar{\vc{v}}(t))|\right>.
\end{align*} 
We use \eqref{imp} and \eqref{support} to conclude 
\begin{align}\label{today0006}
	\int_0^\tau \int_{\mathbb{T}^d}\left<\nu_{t,x};\omega_1^2(s)(s-r(t,x))^2\right>{\rm d}x{\rm d}t&\le \int_0^\tau \int_{\mathbb{T}^d}\left<\nu_{t,x};\omega_1^2(s)(s^2 + r^2(t,x))\right>{\rm d}x{\rm d}t\notag\\&\le\int_0^\tau \int_{\mathbb{T}^d}\left<\nu_{t,x};\omega_1^2(s)(1+s^\gamma + \|r\|_{L^\infty([0,T]\times\mathbb{T}^d)}^2)\right>{\rm d}x{\rm d}t \notag\\&\le C (1+\|r\|_{L^\infty([0,T]\times\mathbb{T}^d)}^2)\int_0^\tau \int_{\mathbb{T}^d}\left<\nu_{t,x}; \omega_1^2(s)(1+s^\gamma)\right>{\rm d}x{\rm d}t\notag\\&\le\,C\,\int_{0}^\tau\mathcal{E}_{mv} \left( s, \vu \ \Big| r , \vc{U} \right)(t){\rm d}t.
\end{align}
We utilize \eqref{compatibility first} to get 
\begin{align}\label{today0007}
&\int_0^\tau\int_{\mathbb{T}^d} \left<\nu_{t,x};|\big(\vU(t)-\bar{\vU}(t)\big)-\big(\vc{v}-\bar{\vc{v}}(t)\big)|^2\right>{\rm d}x{\rm d}s\notag\\&\le C_\lambda \int_0^\tau\int_{{\mathbb{T}^d}\times\mathbb{T}^d}\left<\nu_{t,x,y}^2; ss'|(\vc{v}-\vc{v}')-(\mathbf{U}(t,x)-\mathbf{U}(t,y))|^2 \kappa_\lambda(x-y)\right>\dx\dy\dt + C_\lambda\mathcal{D}(\tau).
\end{align}
We take the help of \eqref{L^2 boundedness}, \eqref{assumption}, \eqref{imp} and \eqref{support} to conclude that 
\begin{align}\label{today0008}
	\int_0^\tau\int_{\mathbb{T}^d} \left<\nu_{t,x},\omega_1(s)|(s-r)||(\bar{\vU}(t)-\bar{\vc{v}}(t))|\right>{\rm d}x{\rm d}t&\le C\,\sup_{t\in[0,T]}|\bar{\vU}(t)-\bar{\vc{v}}(t)|\int_0^\tau\int_{\mathbb{T}^d} \left<\nu_{t,x},\omega_1(s)(1+s^\gamma)\right>{\rm d}x{\rm d}t\notag\\&\le\,C\int_0^\tau \mathcal{E}_{mv} \left( s, \vu \ \Big| r , \vc{U} \right)(t){\rm d}t.
\end{align}
For the term $I_{\omega_2}$, we use \eqref{imp}, \eqref{today0001} and \eqref{support} to obtain
\begin{align}\label{today00009}
	I_{\omega_2}&\le\,C\int_{0}^\tau\int_{\mathbb{T}^d} \left<\nu_{t,x};\omega_2^2(s)\frac{(s-r(t,x))^2}{s}\right>{\rm d}x{\rm d}t + \int_0^\tau \int_{\mathbb{T}^d}\left<\nu_{t,x};s|\vU(t,x)-\vc{v}|^2\right>{\rm d}x{\rm d}t\notag\\&\le
	C\int_0^\tau \mathcal{E}_{mv} \left( s, \vu \ \Big| r , \vc{U} \right)(t){\rm d}t.	
\end{align}
We add \eqref{today0005}-\eqref{today00009} to yield for any $\delta_1>0$
\begin{align}\label{today0009}
	|I_{(1-\psi)}|\le&\,C_\lambda(\delta_1)\,\int_0^\tau \mathcal{E}_{mv} \left( s, \vu \ \Big| r , \vc{U} \right)(t){\rm d}t\notag\\& +  \delta_1 \int_0^\tau\int_{{\mathbb{T}^d}\times\mathbb{T}^d}\left<\nu_{t,x,y}^2; ss'|(\vc{v}-\vc{v}')-(\mathbf{U}(t,x)-\mathbf{U}(t,y))|^2 \kappa_\lambda(x-y)\right>\dx\dy\dt + \delta_1\mathcal{D}(\tau).
\end{align}
Finally, we can use \eqref{today0004} and \eqref{today0009} to get the estimate \eqref{today000}.

\noindent
\textbf{Step 3.}
On the other hand, to estimate the term $I_4$, choose $\delta>0$ small enough such that $[\inf r, \sup r]\subset[\delta, \delta^{-1}]$. Let $\phi_i\in C_c^\infty((0,\infty))$, for $i=1,2$ such that $0\le\phi_i\le 1$ and
\begin{align}\label{imp2}
	\phi_i(s)=1\qquad\text{for}\,\,s\in [\delta/2, 2\delta]\qquad\text{and}\qquad\mbox{supp}\{\phi\}\subset[\delta/4, \delta^{-1}].
\end{align}
We set $\phi(s,s')=\phi_1(s)\phi_2(s')$.
Now we have
\begin{align}
	\begin{split}
	I_4&=\int_0^\tau\int_{\mathbb{T}^d\times\mathbb{T}^d}\left< \nu_{t,x,y}^2; (ss'-r(x)r(y)) {(\vc{U}(t,x) - \vc{U}(t, y))}\kappa_\lambda(x-y)
	\left((\vc{v} - \vc{v}')- (\vc{U}(t,x)-\vc{U}(t,y)) \right)\right>\dx\dy\dt\notag\\
	&=:J_1^\phi + J_2^\phi,
	\end{split}
\end{align}
where
{ \small
\begin{align*}
 J_{1}^{\phi}&= \int_0^\tau\int_{\mathbb{T}^{2d}} \langle \nu_{t,x,y}^2;\phi(s,s') (ss'-r(t,x)r(t,y)){(\vc{U}(t,x) - \vc{U}(t,y))}  \kappa_\lambda(x-y) 
\left((\vc{v} - \vc{v}')- (\vc{U}(t,x)-\vc{U}(t,y)) \right) \rangle \dx\dy\dt \\
J_2^{\phi}&=\int_0^\tau\int_{\mathbb{T}^{2d}}\left< \nu_{t,x,y}^2;(1-\phi(s,s'))(ss'-r(t,x)r(t,y)) {(\vc{U}(t,x) - \vc{U}(t,y))} \kappa_\lambda(x-y)
	\left((\vc{v} - \vc{v}')- (\vc{U}(t,x)-\vc{U}(t,y)) \right)\right>\dx\dy\dt.
\end{align*}}
\textbf{Step 3(a).}
For the term $J_1^\phi$, we have
{ \small
\begin{align}\label{today102}
J_1^\phi&=\int_0^\tau\int_{\mathbb{T}^{2d}}\left< \nu_{t,x,y}^2;\phi(s,s') (ss'-s'r(t,x)) {(\vc{U}(t,x) - \vc{U}(t,y))}\kappa_\lambda(x-y)
	\left((\vc{v} - \vc{v}')- (\vc{U}(t,x)-\vc{U}(t,y)) \right)\right>\dx\dy\dt\notag\\
	&+\int_0^\tau\int_{\mathbb{T}^{2d}}\left< \nu_{t,x,y}^2;\phi(s,s') (s'r(t,x)-r(t,x)r(t,y)){(\vc{U}(t,x) - \vc{U}(t,y))}\kappa_\lambda(x-y)
	\left((\vc{v} - \vc{v}')- (\vc{U}(t,x)-\vc{U}(t,y)) \right)\right>\dx\dy\dt\notag\\
	&=:J_{11}^\phi + J_{12}^\phi.
\end{align}}
We estimate the term $J_{11}^\phi$ as follows: for any $\delta_1'>0$
\begin{align*}
	|J_{11}^\phi|&\le\,C(\delta_1')\,\int_0^\tau\int_{\mathbb{T}^d\times\mathbb{T}^d}\left< \nu_{t,x,y}^2;\phi(s,s')s'(s-r(t, x))^2{|\vc{U}(t,x) - \vc{U}(t,y)|^2}\kappa_\lambda(x-y)
	\right>\dx\dy\dt\\&\qquad + \delta_1'\,\int_0^\tau\int_{\mathbb{T}^d\times\mathbb{T}^d}\left< \nu_{t,x,y}^2;\phi(s,s') s' {|(\vc{v} - \vc{v}')- (\vc{U}(t,x)-\vc{U}(t,y))|^2}\kappa_\lambda(x-y)\right>\dx\dy\dt.
\end{align*}
Making use of the support of $\phi$, we can conclude that
\begin{align*}	
	|J_{11}^\phi|&\le\,C_\lambda(\delta_1')\delta^{-1}\|\mathbf{U}\|_{L^\infty([0,T];W^{1,\infty}(\mathbb{T}^d))}^2\,\int_0^\tau\int_{\mathbb{T}^d}\left<\nu_{t,x};\phi_1(s) |s-r(t,x)|^2\right>\dx\dt\\&\qquad + 4\delta^{-1} \delta_1'\,\int_0^\tau\int_{\mathbb{T}^d\times\mathbb{T}^d}\left< \nu_{t,x,y}^2;\phi(s,s') s's {|(\vc{v} - \vc{v}')- (\vc{U}(t,x)-\vc{U}(t,y))|^2}\kappa_\lambda(x-y)\right>\dx\dy\dt.
\end{align*}
We utilize \eqref{imp}, the compact support of $\phi$, and choose $\delta_1'$ such that $4\delta^{-1} \delta_1' < \delta_1$, to obtain 
\begin{align}\label{today103}
|J_{11}^\phi|&\le\,C(\delta_1)\,\int_0^\tau \mathcal{E}_{mv}(s, \vu \big| r,\mathbf{U})(t){\rm d}t\notag\\&\qquad + \delta_1\,\int_0^\tau\int_{\mathbb{T}^d\times\mathbb{T}^d}\left< \nu_{t,x,y}^2; s's {|(\vc{v} - \vc{v}')- (\vc{U}(t,x)-\vc{U}(t,y))|^2}\kappa_\lambda(x-y)\right>\dx\dy\dt.
\end{align}
We estimate the term $J_{12}^\phi$ as follows: for any $\delta'_1>0$
\begin{align*}
	|J_{12}^\phi|&\le\,C(\delta_1')\|r\|_{L^\infty([0,T]\times\mathbb{T}^d)}\,\int_0^\tau\int_{\mathbb{T}^d\times\mathbb{T}^d}\left< \nu_{t,x,y}^2;\phi(s,s')(s'-r(y,t))^2{|\vc{U}(t,x) - \vc{U}(t,y)|^2}\kappa_\lambda(x-y)
	\right>\dx\dy\dt\\&\qquad + \delta_1'\,\int_0^\tau\int_{\mathbb{T}^d\times\mathbb{T}^d}\left< \nu_{t,x,y}^2;\phi(s,s') {|(\vc{v} - \vc{v}')- (\vc{U}(t,x)-\vc{U}(t,y))|^2}\kappa_\lambda(x-y)\right>\dx\dy\dt.
\end{align*}
Since the support of $\phi$ is contained in $[\delta/4,\delta^{-1}]^2$, we can conclude that
\begin{align*}
|J_{12}^\phi|& \le\,C(\delta_1')\int_0^\tau\int_{\mathbb{T}^d\times\mathbb{T}^d}\left< \nu_{t,x,y}^2;\phi(s,s')(s'-r(t,y))^2{|\vc{U}(t,x) - \vc{U}(t,y)|^2}\kappa_\lambda(x-y)
	\right>\dx\dy\dt\\&\qquad + 16\delta_1'\delta^{-2}\,\int_0^\tau\int_{\mathbb{T}^d\times\mathbb{T}^d}\left< \nu_{t,x,y}^2;\phi(s,s')ss'{|(\vc{v} - \vc{v}')- (\vc{U}(t,x)-\vc{U}(t,y))|^2}\kappa_\lambda(x-y)\right>\dx\dy\dt.
\end{align*}
With the help of \eqref{imp}, the compact support of $\phi$, and choosing $\delta_1'$ such that $16\delta^{-2} \delta_1' < \delta_1$, we get 
\begin{align}\label{today104}
	|J_{12}^\phi|&\le\,C_\lambda(\delta_1)\| \mathbf{U}\|_{L^\infty([0,T];W^{1,\infty}(\mathbb{T}^d))}^2\,\int_0^\tau \mathcal{E}_{mv}(s, \vu \big| r,\mathbf{U})(t){\rm d}t\notag \\&\qquad+  \delta_1\,\int_0^\tau\int_{\mathbb{T}^d\times\mathbb{T}^d}\left< \nu_{t,x,y}^2; s's {|(\vc{v} - \vc{v}')- (\vc{U}(t,x)-\vc{U}(t,y))|^2}\kappa_\lambda(x-y)\right>\dx\dy\dt.
\end{align}
\textbf{Step 3(b).}
For the term $J_2^\phi$, we have for any $\delta_1'>0$,
\begin{align*}
&J_2^\phi=\int_0^\tau\int_{\mathbb{T}^d\times\mathbb{T}^d}\left< \nu_{t,x,y}^2;(1-\phi(s,s'))ss'\bigg(1-\frac{r(t,x)r(t,y)}{ss'}\bigg) 
\left((\vc{v} - \vc{v}')- (\vc{U}(t,x)-\vc{U}(t,y)) \right)\right>\\&\qquad\qquad{(\vc{U}(x) - \vc{U}(t,y))}\kappa_\lambda(x-y)\dx\dy\dt\\&\le\,C(\delta_1')\int_0^\tau\int_{\mathbb{T}^d\times\mathbb{T}^d}\left< \nu_{t,x,y}^2;(1-\phi(s,s'))^2ss'\bigg(1-\frac{r(t,x)r(t,y)}{ss'}\bigg)^2 |\vc{U}(t,x) - \vc{U}(t,y)|^2\kappa_\lambda(x-y)\right>\dx\dy\dt\\&\qquad+\delta_1' \int_0^\tau\int_{\mathbb{T}^d\times\mathbb{T}^d}\left< \nu_{t,x,y}^2;ss' 
|(\vc{v} - \vc{v}')- (\vc{U}(t,x)-\vc{U}(t,y))|^2 \kappa_\lambda(x-y)\right>\dx\dy\dt.
\end{align*}
Let $f\in C((0,\infty)^2)$  such that 
\begin{align*}
	f(s,s')=1 \qquad\forall\, s,s' \ge\,c_\vr\qquad\text{and}\qquad f(s,s')=0\qquad\forall\,s,s'< c_\vr/2.
\end{align*}
By using \eqref{supp1} for the support of $\nu^2_{t,x,y}$, we get
\begin{align*}
	&\int_0^\tau\int_{\mathbb{T}^d\times\mathbb{T}^d}\left< \nu_{t,x,y}^2;(1-\phi(s,s'))^2ss'\bigg(1-\frac{r(t,x)r(t,y)}{ss'}\bigg)^2 |\vc{U}(t,x) - \vc{U}(t,y)|^2\kappa_\lambda(x-y)\right>\dx\dy\dt\\
	&=\int_0^\tau\int_{\mathbb{T}^d\times\mathbb{T}^d}\left< \nu_{t,x,y}^2;(1-\phi(s,s'))^2ss'f(s,s')\bigg(1-\frac{r(t,x)r(t,y)}{ss'}\bigg)^2 |\vc{U}(t,x) - \vc{U}(t,y)|^2\kappa_\lambda(x-y)\right>\dx\dy\dt.
	\end{align*}
Since $f(s,s')\bigg(1-\frac{r(t,x)r(t,y)}{ss'}\bigg)^2$ is bounded on $(0,\infty)^2\times\mathbb{T}^d\times\mathbb{T}^d$. We then obtain 
\begin{align*}
	&\int_0^\tau\int_{\mathbb{T}^d\times\mathbb{T}^d}\left< \nu_{t,x,y}^2;(1-\phi(s,s'))^2ss'\bigg(1-\frac{r(t,x)r(t,y)}{ss'}\bigg)^2 |\vc{U}(t,x) - \vc{U}(t,y)|^2\kappa_\lambda(x-y)\right>\dx\dy\dt\\
	&\le\,C
	\int_0^\tau\int_{\mathbb{T}^d\times\mathbb{T}^d}\left< \nu_{t,x,y}^2;(1-\phi(s,s'))^2ss'|\vc{U}(t,x) - \vc{U}(t,y)|^2\kappa_\lambda(x-y)\right>\dx\dy\dt\\
	&\le\,C\,\int_0^\tau\int_{\mathbb{T}^d\times\mathbb{T}^d}\left< \nu_{t,x,y}^2;(1-\phi(s,s'))^2(s^2 + s'^2)|\vc{U}(t,x) - \vc{U}(t,y)|^2\kappa_\lambda(x-y)\right>\dx\dy\dt\\
	&\le C\,\int_0^\tau\int_{\mathbb{T}^d\times\mathbb{T}^d}\left< \nu_{t,x,y}^2;(1-\phi(s,s'))^2(1+s^\gamma)|\vc{U}(t,x) - \vc{U}(t,y)|^2\kappa_\lambda(x-y)\right>\dx\dy\dt\\&\qquad+C\,\int_0^\tau\int_{\mathbb{T}^d\times\mathbb{T}^d}\left< \nu_{t,x,y}^2;(1-\phi(s,s'))^2(1 + s'^\gamma)|\vc{U}(t,x) - \vc{U}(t,y)|^2\kappa_\lambda(x-y)\right>\dx\dy\dt\\
	&\le\,C_\lambda\| \mathbf{U}\|_{L^\infty([0,T];W^{1,\infty}(\mathbb{T}^d))}^2\int_0^\tau\int_{\mathbb{T}^d\times\mathbb{T}^d}\left< \nu_{t,x,y}^2;(1-\phi(s,s'))^2(1+s^\gamma)\right>\dx\dy\dt\\
	&\le\,C_\lambda\| \mathbf{U}\|_{L^\infty([0,T];W^{1,\infty}(\mathbb{T}^d))}^2\int_{0}^\tau \mathcal{E}_{mv}(s,\vu\big|r, \mathbf{U})(t)\dt,
	\end{align*}
where in the last inequality, \eqref{imp} and \eqref{imp2} are used.
It shows that, by choosing $\delta_1'<\delta_1$, we have
\begin{align}\label{today105}
	J_2^\phi\le &C_\lambda(\delta_1)\int_{0}^\tau \mathcal{E}_{mv}(s,\vu\big|r, \mathbf{U})(t)\dt\notag\\& +\delta_1 \int_0^\tau\int_{\mathbb{T}^d\times\mathbb{T}^d}\left< \nu_{t,x,y}^2;ss' 
	|(\vc{v} - \vc{v}')- (\vc{U}(t,x)-\vc{U}(t,y))|^2 \kappa_\lambda(x-y)\right>\dx\dy\dt.
\end{align}
\textbf{Step 3(c).} We add \eqref{today102}-\eqref{today105} to yield 
\begin{align}\label{today106}
 I_4\le& C_\lambda(\delta_1)\int_{0}^\tau \mathcal{E}_{mv}(s,\vu\big|r, \mathbf{U})(t)\dt\notag\\& +\delta_1 \int_0^\tau\int_{\mathbb{T}^d\times\mathbb{T}^d}\left< \nu_{t,x,y}^2;ss' 
 |(\vc{v} - \vc{v}')- (\vc{U}(t,x)-\vc{U}(t,y))|^2 \kappa_\lambda(x-y)\right>\dx\dy\dt.
\end{align}
\textbf{Step 4.}
By choosing small enough $\delta_1>0$ such that $\delta_1 < 1/8$ and summing up \eqref{today101},\eqref{today000} and \eqref{today106}, we deduce from \eqref{RE100} that
\[
\begin{split}
	&\mathcal{E}_{mv} \left( s, \vu \ \Big| r , \vc{U} \right)(\tau)
	+\frac{7}{8} \mathcal{D}(\tau)+\frac{1}{4}\int_0^\tau\int_{\mathbb{T}^d\times\mathbb{T}^d}\left< \nu_{t,x,y}^2; s's {|(\vc{v} - \vc{v}')- (\vc{U}(t,x)-\vc{U}(t,y))|^2}\kappa_\lambda(x-y)\right>\dx\dy\dt\\
	&\leq \intO{ \left< (\nu_{0})_x;  \frac{1}{2} s |\vc{v} - \vc{U}(0, \cdot) |^2 + P(s) - P'(r(0,\cdot)) (s - r(0,\cdot)) - P(r(0,\cdot)) \right> }\\
	&+ C \left( \int_0^\tau   \mathcal{E}_{mv} \left( s, \vu \ \Big| r , \vc{U} \right)(t)\dt + \int_0^\tau (\chi(t)+\xi(t))\mathcal{D}(t) \ \dt \right).
\end{split}
\]
Thus applying Gronwall's lemma, we conclude that
\begin{equation} \label{RE10}
	\begin{split}
		&\mathcal{E}_{mv}  \left( s, \vu \ \Big| r , \vc{U} \right) (\tau) + \frac{7}{8}\mathcal{D}(\tau) +\frac{1}{4}\int_0^\tau\int_{\mathbb{T}^d\times\mathbb{T}^d}\left< \nu_{t,x,y}^2; s's {|(\vc{v} - \vc{v}')- (\vc{U}(t,x)-\vc{U}(t,y))|^2}\kappa_\lambda(x-y)\right>\dx\dy\dt\\
		&\leq C_T \intO{ \left< (\nu_{0})_x;  \frac{1}{2} s |\vc{v} - \vc{U}(0, \cdot) |^2 + P(s) - P'(r(0,\cdot)) (s - r(0,\cdot)) - P(r(0,\cdot)) \right> },
	\end{split}
\end{equation}
for a.a. $\tau \in [0,T]$. This completes the proof of Theorem \ref{second result}.

\color{black}

\end{document}